%% file: Main_V2.tex
\documentclass[final,3p, times]{elsarticle}
\usepackage{amssymb}
\usepackage{amsmath}
\usepackage{graphicx}	
\usepackage{tikz,pgfplots}
\pgfplotsset{compat=newest}
\usepackage{circuitikz}
\usepackage{color}
\usepackage{url}
\usepackage{booktabs}
\usepackage{wasysym}
\usepackage{xcolor}
\selectcolormodel{gray}

\newcommand{\diag}{\operatorname{diag}}
\newcommand{\norm}[1]{\left\| #1\right\|}
\journal{Electric Power Systems Research}

\usepackage{framed} % Framing content
\usepackage{multicol} % Multiple columns environment
\usepackage{nomencl} % Nomenclature package
\makenomenclature
\renewcommand*\nompreamble{\begin{multicols}{2}}
\renewcommand*\nompostamble{\end{multicols}}
\begin{document}
\begin{frontmatter}
\title{A Convex Approximation for the Tertiary Control of Unbalanced Microgrids.}

\author[1]{Diego-Alejandro Ramirez\corref{cor1}}
\ead{dramirez@xm.com.co}

\author[2]{Alejandro Garc\'es}
\ead{alejandro.garces@utp.edu.co}
% ORCID 0000-0001-6496-0594 

\author[2]{Juan-Jos\'e Mora-Fl\'orez}
\ead{jjmora@utp.edu.co}

\cortext[cor1]{Corresponding author}
\address[1]{XM S.A. E.S.P, Cl. 12 Sur N18-168, Medell\'in, Colombia}
\address[2]{Department of electric power systems engineering. Universidad Tecnol\'ogica de Pereira. AA: 97 - Post Code: 660003 - Pereira, Colombia}

\begin{abstract}
This article presents an optimization model for tertiary control in three-phase unbalanced microgrids. This model considers 24h operation and includes renewable energy sources, energy storage devices, and grid code limitations. Power flow equations are simplified using a recently developed approximation based on Wirtinger's calculus.  The proposed model is evaluated both theoretically and practically.  From the theoretical point of view, the model is suitable for tertiary control since it is convex; hence, global optimum, uniqueness of the solution, and convergence of the interior point method are guaranteed. From the practical point of view, the model is simple enough to be implemented in a small single-board computer with low time calculation. The latter is evaluated by implementing the model in a Raspberry-Pi board with the CIGRE low voltage benchmark; the model is also evaluated in the IEEE 123-nodes test  system for power distribution networks.
\end{abstract}
        
\end{frontmatter}

\section{Introduction}
\subsection{Motivation}

Modern power distribution systems include distributed resources like wind and solar generation as well as storage devices. These components can be grouped, forming a microgrid with potential improvements in efficiency and reliability. However, to achieve these improvements, a microgrid must be complemented with an optimization algorithm called tertiary control, which will be referred to as TC in the rest of the paper \cite{ZIA20181033}.

An optimization model for TC requires to include constraints related to renewable energies and energy storage devices and an accurate representation of the power flow equations. Besides, the model must be tailored for real implementation. Therefore, it is convenient to formulate a model that guarantees global optimum, uniqueness of the solution, and convergence of the algorithms. Convex optimization emerges as a suitable alternative in this context.  Nevertheless, the power flow equations introduce non-convex constraints that require linearization, considering a trade-off between precision and convexity \cite{capitanescu}. On the other hand, the model must be simple enough to be implemented in a practical situation,  with a well-defined time frame for the entire process, and considering time calculations and latency of the communications \cite{7929408}.

\subsection{State of the art}

TC is the highest level in the hierarchy for microgrids operation \cite{DELFINO201875}. It defines the optimal operation point for the active and reactive power in each distributed resource and how much energy the microgrid is willing to trade with the primary grid to satisfy the power balance between the load and power generation \cite{hierarchical_control}.

The problem is closely related to the optimal power flow and the management of energy storage devices. The optimal power flow is a classic problem that has regained importance due to the challenges associated with introducing renewable energy and new advances in convex optimization \cite{ZOHRIZADEH2020391}. A modern approach to this problem is based on conic approximations such as semidefinite programming and second-order cone optimization (see \cite{EES-012} and the references therein for a complete review of this subject).  These relaxations transform the non-convex problem into a convex thereof with theoretical and practical advantages related to global optimal, uniqueness of the solution, and fast convergence rate \cite{YUAN2020106414}.  Despite obvious theoretical advances, these algorithms are far from becoming a proper tertiary control in practical applications \cite{BOBO2021106625}. Computational time may be high, especially in semidefinite approximations, and its implementation may require solvers that allow conic programming \cite{VENZKE2020106480}.  In addition, the algorithms should be tailored to be implemented in low-cost single-board computers and not in a desktop computer.

Another classic approach for the optimal power flow problem is using heuristic algorithms \cite{LI2020117314}.  These algorithms have an extensive history in power systems applications, especially in planning problems.  However, they are not suitable for real-time operation problems.  Computational time is cumbersome in these algorithms, and convergence is not guaranteed \cite{8809848}. Besides, the abuse of heuristic analogies in optimization problems has been criticized in the scientific community \cite{critica_metaheuristica}. Other methods based on artificial intelligence have also been proposed, especially reinforcement learning \cite{9345996}  and neural networks \cite{tertiary3}.  These types of algorithms are promising for real-time implementation, although they may be improved by theoretical analysis of the optimization problem \cite{venzke2020learning}.

TC must consider power flow constraints and energy management of both renewable generation and energy storage resources \cite{GROPPI2021110183}. This aspect is usually studied under the name of energy storage management \cite{tertiary4}. The problem may allow 24h operation but tends to neglect the network effects. The stochastic nature of the problem has been also studied \cite{tertiary5} in both the optimal power flow and the energy management problems. 

\subsection{Contribution}

This paper presents a model based on a linearization recently proposed in \cite{diego_alejandro_general_meeting}. The contribution of the article can be summarized as follows:
\begin{itemize}
    \item An optimization model for 24h operation, considering unbalanced operation, realistic models of the loads and, renewable and storage devices.
    \item A convex approximation of the batteries losses and the exponential models of loads (
    previous linearizations considered constant impedance, current, and power, known as ZIP models and did not consider the non-linear model of the batteries).
    \item A static reserve model that allows the transition from grid-connected to island operation.
    \item An implementation in a small single-board computer (a Raspberry Pi) that demonstrates the proposed algorithm can be implemented in practice at a low cost.
\end{itemize}

Table \ref{tab:estado_arte} shows a comparison with recent references, considering the main aspects for TC in microgrids. $\CIRCLE$ indicates an aspect that was fully considered in the reference, $\RIGHTcircle$ indicates an aspect that was partially considered, and $\Circle$ indicates an aspect that was not considered in said reference.

\begin{table*}[t]
    \centering
    \caption{Aspects considered in recent references about tertiary control for microgrids. }
    \label{tab:estado_arte}
    \begin{tabular}{lccccc}
    \toprule
        Aspect    & Conic relaxations & Heuristics & AI & Energy Management & Proposed\\        
        Reference & \cite{ZOHRIZADEH2020391} to \cite{BOBO2021106625} & \cite{LI2020117314} to \cite{8809848}  & \cite{tertiary3} to \cite{venzke2020learning} & \cite{GROPPI2021110183} to \cite{tertiary5} & approach \\
  \midrule
    Convergence guarantee    & $\CIRCLE$ & $\Circle$ & $\RIGHTcircle$ & $\Circle$      & $\CIRCLE$ \\
    Global optimum           & $\CIRCLE$ & $\RIGHTcircle$ & $\RIGHTcircle$ & $\Circle$      & $\CIRCLE$ \\
    Energy storage           & $\Circle$ & $\Circle$ & $\Circle$      & $\CIRCLE$      & $\CIRCLE$ \\
    Practical oriented       & $\RIGHTcircle$ & $\Circle$ & $\CIRCLE$      & $\RIGHTcircle$ & $\CIRCLE$ \\
    Implementation           & $\Circle$ & $\Circle$ & $\Circle$      & $\Circle$      & $\CIRCLE$ \\
    Stochastic model         & $\Circle$ & $\Circle$ & $\Circle$      & $\CIRCLE$      & $\RIGHTcircle$ \\
    \bottomrule
    \end{tabular}
\end{table*}

\subsection{Organization of the paper}

The paper is organized as follows:  Section 2 shows the architecture of the proposed TC. A general representation of three-phase microgrids and their main components, such as wind turbines, solar panels, and energy storage devices, is presented in Section 3.  Next, a general formulation of the power flow equations for grid-connected operation is presented in Section 4.  The proposed Wirtinger linearization is described in Section 5, resulting in a convex model. The CIGRE microgrid benchmark results are presented in Section 6, followed by its software implementation. Finally, conclusions and recent references are presented.

\section{Tertiary control in microgrids}
\label{sec:tertiary}
Microgrids are usually operated using a hierarchical control that emulates the automatic generation control of conventional power systems. This hierarchy consists of three controls known as primary, secondary, and tertiary control.  Primary control aims to synchronize and stabilize the microgrid, while secondary control aims to achieve nominal frequency; finally, tertiary control aims to optimize the operation.  Therefore, TC operates in a stationary or quasi-stationary state. 

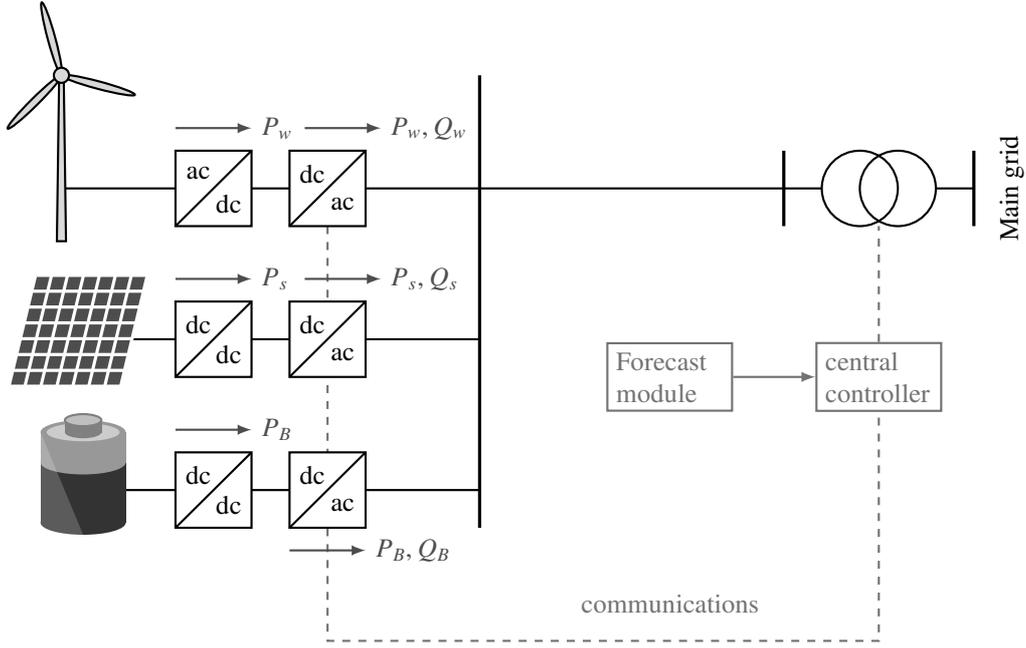
\begin{figure*}[bt]
    \centering
    \input{Figures/Tikz_Convertidores.tex}
    \caption{Integration of wind, solar and energy storage devices through power electronic converters, including TC.}
    \label{fig:power_electronics}
\end{figure*}

A typical architecture for TC is shown in Figure \ref{fig:power_electronics}. Renewable resources and energy storage devices are integrated through power electronics converters to control active and reactive power.  However, this control is limited by the converters' capacity and the primary resource availability (i.e., solar irradiance and wind speed).  Besides, the state of charge of the batteries (SOC) must be managed to minimize power loss and/or reduce costs.  Therefore, an optimization problem requires to be executed by a central controller in real-time.

Real-time may means different things, according to the application.  For a primary control, real-time is an action executed in the order of seconds or even milliseconds, whereas for the tertiary control, real-time is in the order of minutes.  The tertiary control is, therefore, a quasi-dynamic model that works in an entire day.

Different types of communications technologies can be used in this architecture, for example, WiFi or ZigBee
\cite{en12152926}.  In any case, delays or latency must be considered; these delays include those associated with the optimization model and those associated with the communications. The maximum latency experienced by the ZigBee technology for a distance of 50 m with 50-byte packets was 18 ms, according to \cite{latencia2}. Therefore, the primary source of delay is the execution time of the optimization model.  Therefore, fast algorithms with convergence guaranteed are required. 

The central controller can be any device, such as a personal computer or an industrial computer; however, it is advisable an embedded system dedicated explicitly for TC. This system must be small and reliable to be placed in the point of common coupling, that is, the point where the microgrid is electrically connected to the electric distribution system.  The optimization model requires to be accurate but straightforward enough to be implemented in this device.  This is the main objective of the model presented in the next section.

\section{Modelling}

\begin{table*}[htb]
    \begin{framed}
        \printnomenclature
    \end{framed}
\end{table*}

\subsection{Three-phase grid}

\begin{figure}[htb]
    \centering
    \footnotesize
    \input{Figures/tikz_hypergraph.tex}
    \caption{Hypergraph representation for line sections in microgrids}
    \label{fig:hypergraph}
\end{figure}
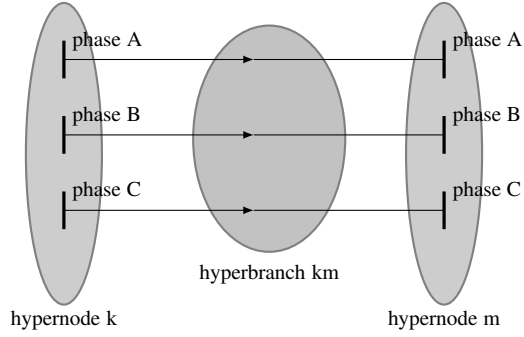

A three-phase microgrid is represented as a connected hypergraph $\mathcal{G} = \left\{\mathcal{H},\mathcal{E} \right\}$, where $\mathcal{H}$ represents the set of hypernodes and $\mathcal{E}\subseteq \mathcal{N}\times\mathcal{N}$ represents the hyperbranches.  Each hypernode and hyperbranch has three components that represent phases $\left\{ A,B,C\right\}$ as depicted in Fig \ref{fig:hypergraph}.
\nomenclature{$\mathcal{H}$}{Set of three-phase nodes (hypernodes)}
\nomenclature{$\mathcal{E}$}{Set of three-phase branches (hyperbranches)}
Along this section and for the rest of the paper, matrices and vectors are represented in capital letter whereas entries of these matrices/vectors are represented in lower case letters. Unless otherwise specified, all variables are defined in the complex domain;  $\overline{x}$ represents the complex conjugate of $x$, and $X^\top$ is the conjugate transpose of $X$. Subscripts represent node and/or time whereas superscripts represent labels of the variables, if they are text, and, exponentiation if they are numeric.

Circuit variables in each hyperbranch $l\in\mathcal{E}$ are represented by a $3\times 3$ admittance matrix as given in \eqref{eq:hyperbrances}.
\begin{equation}
    I_l = Y_l V_l
    \label{eq:hyperbrances}
\end{equation}
These matrices can be grouped together in a block diagonal matrix $Y_\mathcal{E}$ that relates the vector of three-phase voltages with the vector of three-phase currents, as given in \eqref{eq:yprim}.
\begin{equation}
    I_\mathcal{E} = Y_\mathcal{E} V_\mathcal{E}
    \label{eq:yprim}
\end{equation}
where $V_\mathcal{E}$ and $I_\mathcal{E}$ are vectors where the first elements correspond to the phase $A$, next the phase $B$ and finally the phase $C$. An incidence matriz $A\in \mathbb{R}^{\mathcal{H}\times\mathcal{E}}$ is created for the hypergraph, where $a_{kl}=1$ if the hyperbranch $l\in\mathcal{E}$  is in the direction $k\rightarrow m$ and $a_{kl}=-1$ if it is in the direction $m\rightarrow k$. The other entries of the matrix are zero.  This matrix requires to be increased in order to consider three phases in each hypernode. Therefore, the following expressions are obtained:
\begin{align}
    I_\mathcal{H} &= (\mathbb{I}_3\otimes A)I_\mathcal{E} \\
    V_\mathcal{E} &= (\mathbb{I}_3\otimes A)^\top V_\mathcal{H} \\
    I_\mathcal{H} &= (\mathbb{I}_3\otimes A)Y_\mathcal{E}((\mathbb{I}_3\otimes A))^\top V_\mathcal{H} = Y_\mathcal{H} V_\mathcal{H} 
\end{align}
where $\mathbb{I}_3$ is the identity matrix of size $3\times 3$, $Y_\mathcal{H}$ is the three-phase admittance matrix, $\otimes$ is the Kronecker product and $(\cdot)^\top$ represents the transpose conjugate.
\nomenclature{$Y_\mathcal{H}$}{Three-phase admittance matrix}
\nomenclature{$\mathbb{I}_3$}{Identity matrix of size $3\times 3$}
\nomenclature{$I_\mathcal{E}$}{Vector of branch currents}
\nomenclature{$V_\mathcal{E}$}{Vector of branch voltages}
\nomenclature{$I_\mathcal{H}$}{Vector of nodal currents}
\nomenclature{$V_\mathcal{H}$}{Vector of nodal voltages}
\nomenclature{$A$}{Node-branch incidence matrix}
\nomenclature{$\otimes$}{Kronecker product}
\nomenclature{$\mathcal{S}$}{Slack hypernode}
\nomenclature{$\mathcal{N}$}{Hypernodes different from slack, i.e $\mathcal{N}=\mathcal{H}-\mathcal{S}$}

The hypernode set is divided in two new sets $\mathcal{H}=\left\{\mathcal{S},\mathcal{N} \right\}$ where $S$ represents the slack hypernode and $\mathcal{N}$ are the rest of hypernodes.  Therefore, the model of the grid is given by the following matrix equations:
\begin{align}
    I_\mathcal{S} &= Y_\mathcal{SS}V_\mathcal{S} + Y_\mathcal{SN}V_\mathcal{N} \\
    I_\mathcal{N} &= Y_\mathcal{NS}V_S + Y_\mathcal{NN}V_\mathcal{N} 
\end{align}
where
\begin{equation}
    Y_\mathcal{H} = \left(\begin{array}{cc} Y_\mathcal{SS} & Y_\mathcal{SN} \\ Y_\mathcal{NS} & Y_\mathcal{NN}\end{array} \right)
\end{equation}
Notice that $\mathcal{S}$ is size three since there are three slack nodes in the system, corresponding to each phase.  The voltage $V_\mathcal{S}$ is given by \eqref{eq:vslackk}:
\begin{equation}
    V_\mathcal{S} = v^\text{nom}\left(\begin{array}{c}1 \\ e^{-2\pi/3j} \\ e^{2\pi/3j}  \end{array}\right)
    \label{eq:vslackk}
\end{equation}
where $v^\text{nom}$ is the nominal voltage in per unit, measured in the point of point of common coupling. On the other hand, voltages $V_\mathcal{N}$ are variables in the optimization model and thus, the current in each node of the grid is in function of these voltages
\begin{equation}
    i_k = \sum_{s\in \mathcal{S}} y_{sk}v_s + \sum_{m\in \mathcal{N}} y_{km}v_m, \; \forall k \in \mathcal{N}
\end{equation}
multiplying by $\overline{v}_k$ an expression for the nodal powers is obtained as follows:
\begin{equation}
    \overline{s}_{kt} = \sum_{m\in\mathcal{S}}\overline{v}_{kt}y_{mk}v_{mt} + \sum_{m\in\mathcal{N}} \overline{v}_{kt}y_{km}v_{mt}, \; \forall k \in \mathcal{N}
    \label{eq:ecuaciones_flujo_de_carga_completas}
\end{equation}
where $\overline{v}_k$ represents the complex conjugate of $v_k$ and $s_k=p_k+jq_k$ is the total apparent power in node $k$. A subscripts $t$ was added to all variables to indicate the time since the model is formulated for one day operation.

Each node may have different power injections such as wind and solar generation, energy storage and, loads.  Therefore, the following balance of power in each node is obtained:
\begin{equation}
s_{kt} = s^\text{wind}_{kt}+s^\text{solar}_{kt}+s^\text{battery}_{kt}- s^\text{load}_{kt}
\label{eq:balance_potencia_total}
\end{equation}
Total losses of the grid are calculated as given in \eqref{eq:Pl}.
\begin{equation}
    p^\text{loss}_t = \operatorname{real}\left( V_{\mathcal{S}t}^\top Y_\mathcal{SS} V_{\mathcal{S}t}+2V_{\mathcal{S}t}^\top Y_\mathcal{SN} V_{\mathcal{N}t} + V_{\mathcal{N}t}^\top Y_\mathcal{NN} V_{\mathcal{N}t}\right)
    \label{eq:Pl}
\end{equation}
\nomenclature{$s_{kt}$}{Apparent power in node $k$ at time $t$}
\nomenclature{$p^\text{loss}_t$}{Total power loss at time $t$}
\nomenclature{$p^\text{grid}_t$}{Total power supplied for the main grid at time $t$}
Further, the power supplied by the main grid to the microgrid is defined as follows
\begin{equation}
    p^\text{grid}_t=\operatorname{real}\left(\left(Y_\mathcal{SS} V_{\mathcal{S}t} + Y_\mathcal{SN}V_{\mathcal{N}t}\right)^\top V_{\mathcal{S}t}\right)
    \label{eq:potencia_grid}
\end{equation}

This variable defines the interchange of power with the main grid, and hence, it is essential for the optimization model.

\subsection{Three-phase exponential load model}

Three-phase loads are represented by a general model that considers constant power, constant current, and constant impedance \cite{ZIP}. This model can be described in terms of a constant $\alpha\in\left\lbrace0,1,2\right\rbrace$ as given in \eqref{eq:modelo_carga} where $0$ corresponds to constant power, $1$ to constant current and $2$ to constant impedance.  Fractional values of $\alpha_k$ are also allowed in the model.
\begin{equation}
    s^\text{load}_{kt} = s^\text{ZIP}_{kt} \left\|\frac{v_{kt}}{v^\text{nom}}\right\|^{\alpha_k}
    \label{eq:modelo_carga}
\end{equation}
\nomenclature{$\alpha_k$}{Exponential model of load $k$}
\nomenclature{$v^\text{nom}$}{Nominal phase-to-neutral voltage}
A time series for $s^\text{load}_k$ is required in order to obtain the 24h operation.
\label{sec:load_model}

\subsection{Photovoltaic units}
The model that represents the power supplied by the photovoltaic units (PVs) during a period is presented in \eqref{eq:solar}
\begin{equation}
    p_{kt}^\text{solar} \leq \rho_{k} \ \psi_t
    \label{eq:solar}
\end{equation}
\nomenclature{$\psi_t$}{Irradiance at time $t$}
\nomenclature{$\rho_{k}$}{Productivity coefficient associated with the solar panel $k$}
\nomenclature{$s_{kt}^\text{solar}$}{Power generated by photovoltaic unit $k$ at time $t$}
where $p_{kt}^\text{solar}$ power that the photovoltaic generator can inject in the period $t$,  $\rho_k$ is the productivity coefficient associated to the solar panel, and $\psi_t$ corresponds to the irradiance measured in $(\text{W}/\text{m}^2)$, perpendicular to the panel for the period $t$. The model requires a time series of $\psi$ for 24 h.  Power electronic converters are allow to reduce their production if required by the optimization model there in the inequality.

\subsection{Wind turbine model}
Generally, a model for the generated wind power may be obtained using a wind turbine profile obtained experimentally. In this case, the proposed model represents an approximation of the active power supplied by the wind turbine; thus the wind power may be calculated based on the wind speed and the wind turbine power coefficient as follows:
\begin{equation}
    p_{kt}^\text{wind}\leq\left\{\begin{array}{ccc}
    p^\text{nom}_k\left(\frac{w_t }{{w^\text{nom}_k}}\right)^3, & \text{if} & 0\leq w_t \leq w^\text{nom}_k\\
    \\ p^\text{nom}_k, & \text{if} & w^\text{nom}_k \leq w_t \leq w^\text{max}_k\\
    \\ 0, & \text{if} & w^\text{max}_k \leq w_t
    \end{array}
\right\}
\label{eq:wind_generation}
\end{equation}
\nomenclature{$s_{kt}^\text{wind}$}{Power generated by wind turbine $k$ at time $t$}
where $p_{kt}^\text{wind}$ corresponds to the maximum power that the wind turbine can deliver in the period of time $t$, the wind speed is a time variable described by $w_t$, and $p^\text{nom}$ is the rated power that each of the wind turbines can inject to the system. $w^\text{nom}$ represents the rated wind speed at which the turbines generate its rated power, and the maximum wind speed that the turbines tolerates is represented by $w^\text{max}$. The model requires a time series of the wind velocity during a period of time.
\nomenclature{$w_t$}{wind speed at time $t$}
\subsection{Battery energy storage}

The following equation proposes a simplified dynamic model of the battery energy storage $k$ at each time $t$.

\begin{equation}
    e_{kt} = e_{kt-1}+\left(p^\text{char}_{kt}-p^\text{disch}_{kt}\right) \Delta t
    \label{eq:energy_balance_battery}
\end{equation}
\nomenclature{$e_{kt}$}{Energy stored in the battery $k$ at time $t$}
where $p^\text{char}_{kt}$ is the charging power, $p^\text{disch}_{kt}$ is the power at discharge, and $e_{kt}$ is the energy stored by the batteries, considering a positive value if it injects power to the microgrid and a negative value for the charging mode. $\Delta t$ is the time discretization (usually 1h). A quadratic model of the losses in the battery and the converter is proposed, considering efficiencies of charge and discharge as follows:
\begin{align}
    n^\text{char}_{kt} &= a^\text{char}_k \left(p^\text{char}_{kt}\right)^2 + b^\text{char}_kp^\text{char}_{kt} + c^\text{char}_k \label{eq:bateria_carga} \\
    n^\text{disch}_{kt} &= a^\text{disch}_k \left(p^\text{disch}_{kt}\right)^2 + b^\text{disch}_kp^\text{disch}_{kt} + c^\text{disch}_k \label{eq:bateria_descarga}
\end{align}
this model is taken from data-sheets of batteries (see for example \cite{data_sheet}).  It is important to remark that $a^\text{charge}$ and $a^\text{discharge}$ are always positive.
Then, the active power delivered by the battery energy storage is given by the following equations:
\begin{equation}
    p^\text{char}_{kt} - p^\text{disch}_{kt} = \operatorname{real}\left(s^\text{battery}\right) 
    \label{eq:cargadescarga}
\end{equation}

\subsection{Power and voltage limits}

Photovoltaic units, wind turbines, and battery energy storage devices are integrated into the microgrid through power electronic converters, as depicted in Fig
\ref{fig:power_electronics}.  These converters can generate or consume reactive power on the assumption that the current and apparent power is below its maximum.

Therefore, the apparent power is limited by the following constrain:
\begin{equation}
\left\| s_k^\text{device} \right\| \leq s_{\max}^\text{device}, \; \forall \text{device}\in\left\{\text{wind,solar,battery} \right\} 
\label{eq:power_electronics}
\end{equation}
Notice that this constraint is convex.

On the other hand, voltage limits are considered by the following constraint:
\begin{equation}
    \left\| v_k-v^\text{nom} \right\| \leq \delta v^\text{nom} 
    \label{eq:voltage_limits}
\end{equation}
where $\delta$ is the maximum voltage deviation according to the grid code.

\subsection{Static reserve}

A \textit{static reserve} is proposed in the model. This new concept is strongly related to the energy stored by the batteries. It is applied because it needs to keep the system secure and in equilibrium when the system changes from grid-connected to island operation.  The static reserve guarantees a time of safe operation in which the demand of the microgrid can be supplied by the available generation and energy stored by the batteries. Taking this criterion into account, a deficit $d_t$ of power provided by generators is defined as a function of the demand and the power injected by the slack node to the microgrid as given in \eqref{eq:deficitt},

\begin{equation}
    d_t = \sum_{k\in\mathcal{N}} s^\text{load}_k-p^\text{grid}_t
\label{eq:deficitt}    
\end{equation}

Therefore, the energy stored by the batteries is limited by the following constraints:

\begin{equation}
    \sum_{k\in\mathcal{N}^\text{battery}}e_{kt}\geq \operatorname{real}(d_t)
    \label{eq:deficitt2}
\end{equation}
Where $\mathcal{N}^\text{battery}\subseteq\mathcal{N}$ is the subset of nodes that have battery energy storage.

\section{Optimization problem}

The proposed power flow requires a time series for wind speed, solar irradiance, and power demand. The optimization model search for optimal use of energy storage devices as given below, where all variables depend on the time, and therefore, the subindex $t$ is omitted in most of the equations. 

The model seeks to minimize the costs $c_t$ of the energy supplied by the primary grid. The model of each component and the grid itself are also considered. For the sake of completeness, the entire model is presented in Table \ref{tab:non_convex_model} including the type of equation (notice that quadratic equality constraints are also non-convex).

\begin{table}[hbt]
    \centering
    \caption{Non-convex model for the tertiary control  (NC-TC Model)}
    \label{tab:non_convex_model}
    \begin{tabular}{ccc}
    \toprule
    Equation & Type & Characteristic  \\
    \hline 
    $\min \; \sum_t c_tp^\text{grid}_t$     & objective function  & affine \\
    \hline      
    \eqref{eq:ecuaciones_flujo_de_carga_completas}  & load flow & non-convex         \\
    \eqref{eq:balance_potencia_total} & power balance & affine \\
    \eqref{eq:Pl}  & power loss & quadratic equality \\
    \eqref{eq:potencia_grid} & grid power & quadratic equality \\
    \eqref{eq:modelo_carga} & load model & non-convex \\
    \eqref{eq:solar} & solar generation & affine \\
    \eqref{eq:wind_generation} & wind generation & affine \\
    \eqref{eq:energy_balance_battery} & energy balance & affine \\
    \eqref{eq:bateria_carga} & battery charging & quadratic equality \\
    \eqref{eq:bateria_descarga} & battery discharge & quadratic equality \\
    \eqref{eq:cargadescarga} & battery power & affine \\
    \eqref{eq:power_electronics} & converter capacity & second-order cone \\
    \eqref{eq:voltage_limits} & voltage limits & second-order cone \\
    \eqref{eq:deficitt} & power deficit & affine \\
    \eqref{eq:deficitt2} & static reserve & affine \\
    \bottomrule 
    \end{tabular}
\end{table}

\nomenclature{$n^\text{char}_{kt}$}{Losses in the energy storage during charging}
\nomenclature{$n^\text{disch}_{kt}$}{Losses in the energy storage during discharging}
\nomenclature{$p^\text{char}_{kt}$}{Power in the battery during charging}
\nomenclature{$p^\text{disch}_{kt}$}{Power in the battery during discharging}
\nomenclature{$e_{kt}$}{Energy in the battery $k$ at time $t$}

The model presented in Table \ref{tab:non_convex_model} is called, in this paper, as non-convex tertiary control model  (hereafter called NC-TC Model).  This model contains some complex variables (e.g $s_k$,$v_k$ and $v_m$) which simplifies its representation.  However, it is important to remark that this is only a representation since an optimization model requires to be defined in an ordered set (for example, the set $\mathbb{R}^n$).  Therefore, each equality constraint with complex variables must be separated into real and imaginary parts.  This complex representation allows a simple formulation of the linearizations using Wirtinger's calculus.  The resulting linearizations can be written directly in complex form, in packages such as \textit{cvxpy} as explained in the next section.

\subsection{Convex identification}

NC-TC model is non-linear and non-convex, and therefore, a convex approximation is required. However, some parts of the model are already convex, for instance, affine equations.  Power loss, grid power, and loss in the battery energy storage devices are all quadratic equality constraints.  These equations define, in principle, non-convex sets.  However, they can be transformed into convex sets by replacing equality for inequality constraints and noticing that these quadratic forms are all convex. This approximation is suitable since the objective function includes minimizing $p^\text{grid}$ and hence, the inequality gap tends to be minimized. 

The only remaining non-convex constraints are \eqref{eq:ecuaciones_flujo_de_carga_completas} and \eqref{eq:modelo_carga}. These constraints will be linearized by using Wirtinger calculus as presented in the next subsection.

\section{Convex model}
\subsection{Wirtinger's linearization}

In general, an equality constraint is non-convex unless it is an affine equation. Therefore, a linearization is recommended in order to approximate the model to a convex set. 

Let us consider a complex variable $z=x+jy$ and a complex function $f(z)=u+jv$. The Wirtinger's derivate and the conjugate Wirtinger's derivate are defined as follows:
\begin{align}
    \frac{\widehat{\partial}f}{\partial z}&=\frac{1}{2}\left(\frac{\partial u}{\partial x}+\frac{\partial v}{\partial y} \right)+\frac{j}{2}\left(\frac{\partial v}{\partial x}-\frac{\partial u}{\partial y} \right)\\
    \frac{\widehat{\partial} f}{\partial \overline{z}}&=\frac{1}{2}\left(\frac{\partial u}{\partial x}-\frac{\partial v}{\partial y} \right)+\frac{j}{2}\left(\frac{\partial v}{\partial x}+\frac{\partial u}{\partial y} \right)
\end{align}
\nomenclature{$\partial$}{Conventional derivative}
\nomenclature{$\widehat{\partial}$}{Wirtinger's derivative}
These operators are very similar to the conventional derivatives, despite the fact that they are not derivatives in the Cauchy-Riemann sense (see \cite{Wirtinger} and \cite{Wirtinger2} for more details).  More importantly, these operators allow a linearization on the complex numbers for non-holomorphic functions such as \eqref{eq:ecuaciones_flujo_de_carga_completas} and \eqref{eq:modelo_carga}.

A complex function linearization in terms of the Wirtinger's operators is defined as follows:
\begin{equation}
    f \approx  f(z_0) + \frac{\widehat{\partial} f}{\partial z} \Delta z + \frac{\widehat{\partial} f}{\partial \overline{z}} \Delta \overline{z}
    \label{eq:wlinearization}
\end{equation}
For example, a function $f=\overline{v}_kv_m$ can be linearized around a point $v_{k0},v_{m0}$ as follows:
\begin{align}
f(\overline{v}_{k},v_{m})&=\overline{v}_{k}v_{m} \\
        &\approx \overline{v}_{k0}v_{m0}+v_{m0}\Delta \overline{v}_{k}+\overline{v}_{k0}\Delta v_{m} \\
        &= \overline{v}_{k0}v_{m0}+v_{m0}(\overline{v}_{k}-\overline{v}_{k0})+\nonumber\\&\ \ \ \ \overline{v}_{k0}(v_{m}-v_{m0}) \\
        &= v_{m0}\overline{v}_{k}+\overline{v}_{k0}v_{m}-\overline{v}_{k0}v_{m0}\label{eq:linearization}
\end{align}
Notice that Wirtinger's derivatives fulfill the basic properties for differentiation known from real-valued analysis concerning the linearity, product rule and composition of functions.

\subsection{Wirtinger's Linearization of the power flow equations}

The linearization presented in \eqref{eq:linearization} can be used for linearizing \eqref{eq:ecuaciones_flujo_de_carga_completas}, as follows:

\begin{equation}
\overline{s}_{kt} = \sum_{m\in\mathcal{S}}\overline{v}_{kt}y_{mk}v_{mt} + \sum_{m\in\mathcal{N}} y_{km}\left(v_{m0}\overline{v}_{k}+\overline{v}_{k0}v_{m}-\overline{v}_{k0}v_{m0}\right)
\label{eq:lineal}
\end{equation}

After simple algebraic manipulations the following matrix representation is obtained:
\begin{eqnarray}
    \overline{S}_{\mathcal{N}t}=\diag(Y_{\mathcal{S}}\cdot V_{\mathcal{S}t})\cdot \overline{V}_{\mathcal{N}t}+\diag(Y_{\mathcal{N}}\cdot V_{\mathcal{N}0})\cdot \overline{V}_{\mathcal{N}t}+ \nonumber\\
    \left(\diag(\overline{V}_{\mathcal{N}0}\right)\cdot Y_{\mathcal{N}})\cdot V_{\mathcal{N}t}-\diag(V_{\mathcal{N}0})\cdot(Y_{\mathcal{N}}\cdot \overline{V}_{\mathcal{N}0}) 
    \label{eq:matrixform}
\end{eqnarray}
Equation \eqref{eq:matrixform} defines an affine space which constitutes an approximation of the power flow. The following matrices are defined in order to simplify the nomenclature:
\begin{align}
    K&=\diag(Y_{\mathcal{S}} V_{\mathcal{S}})+\diag(Y_{\mathcal{N}}\cdot V_{\mathcal{N}0})\\
    L&=\diag(\overline{V}_{\mathcal{N}0})\cdot Y_{\mathcal{N}}\\
    U&=-\diag(V_{\mathcal{N}0})\cdot(Y_{\mathcal{N}}\cdot \overline{V}_{\mathcal{N}0})
\end{align}
Therefore, Constraint \eqref{eq:ecuaciones_flujo_de_carga_completas} is transformed into the following affine equation:
\begin{equation}
    \overline{S}_{\mathcal{N}t}=K\cdot \overline{V}_{\mathcal{N}t}+L\cdot V_{\mathcal{N}t}+U
    \label{eq:simplify}
\end{equation}
Precision of this approximation was evaluated in \cite{diego_alejandro_general_meeting}, where results demonstrated errors less than $1\times 10^{-3}$ in voltages.  It is important to remark that the linearization is defined in the complex domain, therefore, $v_{k0}$ must consider the phase (i.e., $e^{0j}$ for phase A, $e^{-2\pi/3j}$ for phase B and $e^{2\pi/3j}$ for phase C).

\subsection{Wirtinger's linearization for exponential loads}

Equation \eqref{eq:modelo_carga} can be also linearized using Wirtinger's calculus. The non-linear term in this expression is given by \eqref{eq:general3}.
\begin{equation}
    \norm{v_k}^\alpha = (v_k \overline{v}_k)^{\frac{\alpha}{2}}
    \label{eq:general3}
\end{equation}
Linearizing around $v_{k0}$ the following expression is obtained 
\begin{equation}
\begin{split}
               (v_k  \overline{v}_k)^{\frac{\alpha}{2}} &\approx (v_{k0}  \overline{v}_{k0})^{\frac{\alpha}{2}} + \frac{\widehat{\partial} \norm{v_k}^\alpha}{\partial v_k} \Delta v_k + \frac{\widehat{\partial} \norm{v_k}^\alpha}{\partial \overline{v}_k} \Delta \overline{v}_k\\
                &= (v_{k0} v_{k0}^*)^{\frac{\alpha}{2}} + \frac{\alpha}{2}(v_{k0} \overline{v}_{k0})^{\frac{\alpha}{2}-1} (\overline{v}_{k0} \Delta v_k+v_{k0} \Delta \overline{v}_k)
\end{split}
\label{eq:lZIP}
\end{equation}
Considering that the voltages of the system are near of the nominal value, a good point of linearization is $v_{k0}=v^\text{nom}e^{j\phi}$ as the complex nominal voltage, where $\phi$ is $0,-2\pi/3$ or $2\pi/3$ according to the phase. Therefore, the nonlinearity of the ZIP load model is represented linearly by the following expression:
\begin{equation}
\begin{split}
    \left(v_k \overline{v}_k\right)^{\frac{\alpha}{2}}\approx\left(v^\text{nom}  \right)^{\alpha} + \\ 
    \frac{\alpha}{2}\left(v^\text{nom}  \right)^{\alpha-2} \left(v^\text{nom}e^{-j\phi} \left(v_k-v^\text{nom}e^{j\phi}\right)+
    v^\text{nom}e^{j\phi} \left(\overline{v}_k-v^\text{nom}e^{-j\phi} \right)\right)
\end{split}        
\label{eq:linVnom}
\end{equation}
Then, simplifying \eqref{eq:linVnom}, the following equation is obtained:
\begin{align}
    \left(v_k  \overline{v}_k\right)^{\frac{\alpha}{2}}&=\left(v^\text{nom}  \right)^{\alpha} + \nonumber \\ & \frac{\alpha}{2}\left(v^\text{nom} \right)^{\alpha-2} \left(v^\text{nom}e^{-j\phi} v_k+v^\text{nom}e^{j\phi} \overline{v}_k-2\left(v^\text{nom}\right)^2\right)
    \label{eq:ZIPlinearized}
\end{align}
Therefore, \eqref{eq:modelo_carga} is replaced by the following affine equation:
\begin{equation}
S^\text{load}_{\mathcal{N}t} = S^\text{ZIP}_{\mathcal{N}t} \circ \left(M+HV_{\mathcal{N}t} +T\overline{V}_{\mathcal{N}t}\right)
    \label{eq:linealizacion_zip}
\end{equation}
\nomenclature{$\circ$}{Hadamard product}
where $\circ$ represents the Hadamard product and $M,H,T$ are constant vectors given by
\begin{align}
    M &= \operatorname{diag}\left(1-\alpha\right)  \\
    H &= \operatorname{diag}\left(\frac{\alpha}{2 v^\text{nom} e^{j\phi}}\right) \\
    T &= \operatorname{diag}\left(\frac{\alpha}{2 v^\text{nom} e^{-j\phi} }\right)
\end{align}
For the sake of completeness, the model with all the aforementioned approximations is presented in Table \ref{tab:modelo_convexo}.

\begin{table}[hbt]
    \centering
    \caption{Convex model for the tertiary control (Convex-TC)}
    \label{tab:modelo_convexo}
    \begin{tabular}{ccc}
    \toprule
    Equation & Type & Modification  \\
    \hline 
    $\min \; \sum_t c_tp^\text{grid}_t$     & objective function  & affine \\
    \hline      
    \eqref{eq:simplify}  & load flow & non-convex         \\
    \eqref{eq:balance_potencia_total} & power balance & affine \\
    \eqref{eq:Pl}  & power loss & quadratic inequality \\
    \eqref{eq:potencia_grid} & grid power & quadratic inequality \\
    \eqref{eq:linealizacion_zip} & load model & affine \\
    \eqref{eq:solar} & solar generation & affine \\
    \eqref{eq:wind_generation} & wind generation & affine \\
    \eqref{eq:energy_balance_battery} & energy balance & affine \\
    \eqref{eq:bateria_carga} & battery charging & quadratic inequality \\
    \eqref{eq:bateria_descarga} & battery discharge & quadratic inequality \\
    \eqref{eq:cargadescarga} & battery power & affine \\
    \eqref{eq:power_electronics} & converter capacity & second-order cone \\
    \eqref{eq:voltage_limits} & voltage limits & second-order cone \\
    \eqref{eq:deficitt} & power deficit & affine \\
    \eqref{eq:deficitt2} & static reserve & affine \\
    \bottomrule 
    \end{tabular}
\end{table}

This model will be called Convex-TC in the following sections in order to differentiate from the NC-TC model.
Notice that, the objective function includes a quadratic function with a positive definite matrix. This property is guaranteed if the hypergraph is connected. Therefore, the function is strongly convex, and this implies a global optimum and unique solution.  In addition, the quadratic inequality constraints are also convex.

\subsection{Operation under surplus energy limitation}

Some grid codes prevent the sale of surpluses of energy.  In these cases, the model must be limited with an additional constrain in order to prevent the power in the slack node becomes negative.  The optimization model is modified in order to include this constraint as follows: 
\begin{equation}
\begin{split}
    \min \; & \sum_t c_tp^\text{grid}_t \\
    & \operatorname{real}\left(p^\text{grid}_t\right) \geq 0 \\
    & \operatorname{imag}\left(p^\text{grid}_t\right) \geq 0 \\
    & \text{+ all constrains of Convex-TC} 
\end{split}
\end{equation}

When the available energy is higher than the total load, the microgrid remains connected to the main grid, but the active and reactive power interchanged between the grid, and the microgrid is zero.  

\section{Results}\label{sec:Results}
\subsection{CIGRE low voltage benchmark test system}

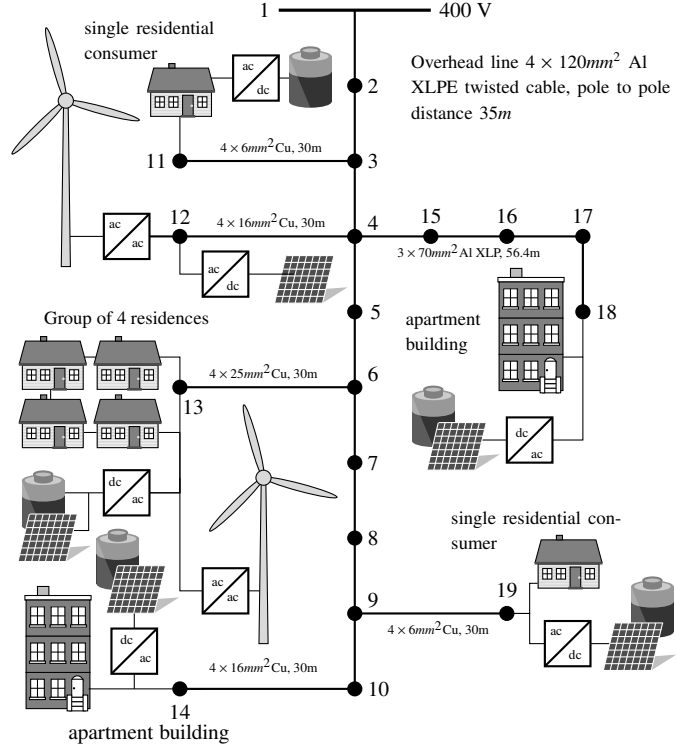
\begin{figure}[bth]
    \centering
    \input{Figures/tikz_cigre.tex}
    \caption{The CIGRE low voltage benchmark test system.}
    \label{fig:cigre}
\end{figure}

A modified version of the benchmark test system for low-voltage microgrids proposed by CIGRE in \cite{CIGRE} was initially used to illustrate the use of the model and its performance. This system is a 19-nodes, typical residential network with a peak power demand of 186.9 kW and a nominal voltage of 400V. The system was modified to include renewable generation and energy storage devices, as depicted in Fig \ref{fig:cigre}.

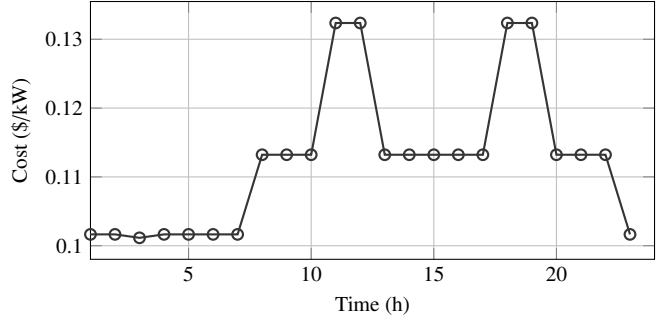
\begin{figure}[bth]
    \footnotesize
		\centering
    \input{Figures/Costo.tex}
    \caption{Prices from XM S.A.S, Colombian national interconnected system operator. Real time pricing for 24 consecutive hours in 2019.}
    \label{fig:Costo}
\end{figure}

The exponential model described in subsection \ref{sec:load_model} was used for the loads, with values of $\alpha$ given in Table \ref{table:load}. The system was analyzed for operation in 24h with intervals of 1h; a variable energy cost $c_t$ was considered as shown in Figure \ref{fig:Costo}.

\begin{table}[bth]
\centering
\caption{Loads description according to the exponential load model}
\label{table:load}
\begin{tabular}{ccc}
\toprule
Node & Peak power demand (W) & $\alpha$ \\
\hline
     11 & 13400 & 2 \\
     13 & 47000 & 0 \\
     14 & 40000 & 2 \\
     18 & 70900 & 0 \\
     19 & 15600 & 1 \\
\bottomrule
\end{tabular}
\end{table}

Two main cases were analyzed, namely:
\begin{description}
\item[Case 1 (surplus sale):] in this case, the slack node can inject and receive power from the microgrid; the optimization model decides the purchase or sale of energy, according to market prices and the generation and demand conditions of the microgrid.

\item[Case 2 (without surplus sale):] in this case, the slack node is limited to only supply power, i.e., the microgrid cannot inject power to the main grid. This type of limitation is due to some grid codes that are restrictive in terms of surplus sales.  
\end{description}

\subsection{Case 1 (surplus sale)} 
\label{subsec:Case1}
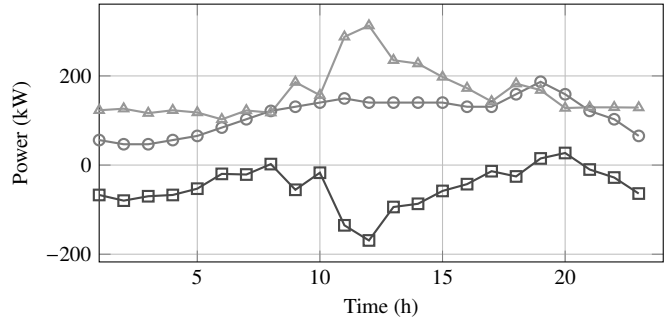
\begin{figure}
    \centering
    \footnotesize
    \input{Figures/GenvsDem_PowSale.tex}
    \caption{Tertiary control the slack node (\ref{tikz:slack}), load demand (\ref{tikz:demand}) and power $P_{ER}$ (\ref{tikz:per}) for the Case 1.}
    \label{fig:GenvsDem_PowSale}
\end{figure}

The behavior of the TC for this case is shown in Fig \ref{fig:GenvsDem_PowSale}, where $p^\text{grid}$ is the power taken or sold to the main grid, and $p^{ER}$ represents the power provided by the photovoltaic units, wind turbines, and batteries. As a result, the microgrid takes power from the main grid in periods where the renewable generation and batteries can not satisfy the load or when the energy has a low price. Eventually, the microgrid sold power to the system when the generator can supply the load, and there is an excess of generation.

\begin{figure}
    \centering
    \footnotesize
    \input{Figures/SvSmax_PowSale.tex}
    \caption{Optimal solar generation (\ref{tikz:ppv}) and available solar generation (\ref{tikz:pvmax})  for the Case 1.  The algorithm always takes the maximum solar generation.}
    \label{fig:SvSmax_PowSale}
\end{figure}
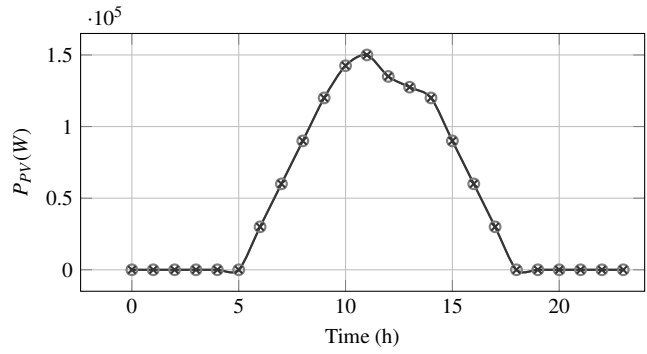

The photovoltaic generation uses all its available resources, injecting the maximum amount of power that the irradiance profile allows in all the periods. The behavior of the photovoltaics units is represented in Fig \ref{fig:SvSmax_PowSale}.

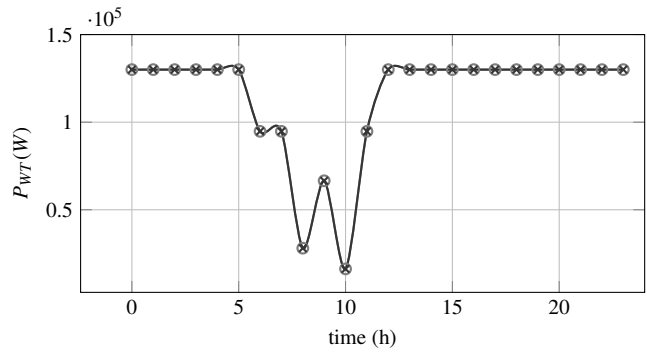
\begin{figure}
    \centering
    \footnotesize
    \input{Figures/WvWmax_PowSale.tex}
    \caption{Optimal wind generation (\ref{tikz:pw}) and available wind generation (\ref{tikz:pwmax}) for Case 1.  The algorithm always take the maximum available power.}
    \label{fig:WvWmax_PowSale}
\end{figure}

Wind turbines have the same behavior as photovoltaic sources, as shown in Fig \ref{fig:WvWmax_PowSale}. These use the available wind speed profile, giving a result where the turbine can provide the maximum power in each period according to the wind speed forecast.

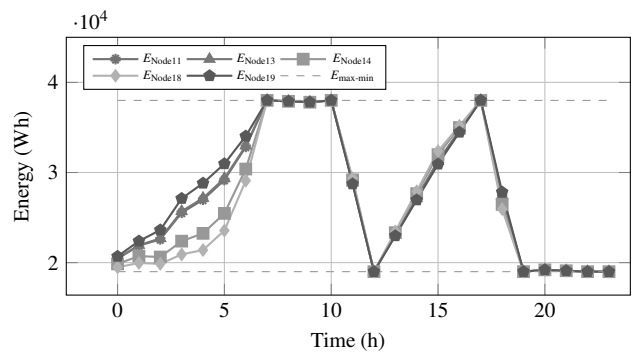
\begin{figure}
    \centering
    \footnotesize
    \input{Figures/Baterias_PowSale.tex}
    \caption{Energy stored in each battery for Case 1, considering the SOC of the batteries between $30\%$ and $100\%$ of the batteries energy.}
    \label{fig:Baterias_PowSale}
\end{figure}

Finally, and considering an initial batteries SOC of $50\%$, these have a response to the load, as shown in Fig \ref{fig:Baterias_PowSale}.

Notice that all batteries have similar behavior; they charge completely in periods of low prices, low demand and/or high generation, to discharge in periods with high demand when the prices are higher. 

\subsection{Case 2 (without surplus sale)} 

The power provided by the slack node is limited in all the periods of the tertiary control.  This constraint is represented as follows:
\begin{equation}
   p^\text{grid}_t \geq 0
\end{equation}

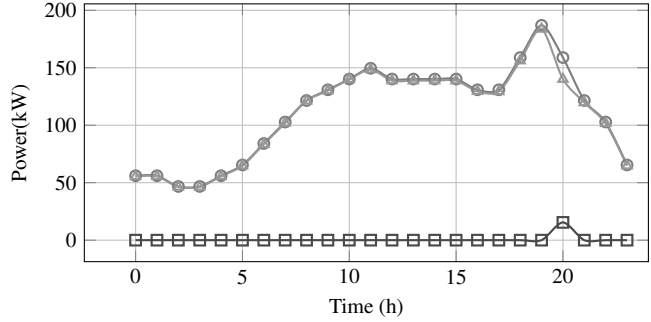
\begin{figure}
    \centering
    \footnotesize
    \input{Figures/GenvsDem.tex}
    \caption{Energy exchange with the main grid (\ref{tikz:pslack}), total generation (\ref{tikz:generation}) and total demand (\ref{tikz:totaldemand}) for Case 2.}
    \label{fig:GenvsDem}
\end{figure}

Therefore, renewable generation never overtakes the load demand, and the microgrid can not inject power into the main grid. In the periods that the demand exceeds the maximum renewable generation and the power that the batteries can provide, the main grid injects power to the microgrid to supply the loads, as shown in Fig \ref{fig:GenvsDem}.

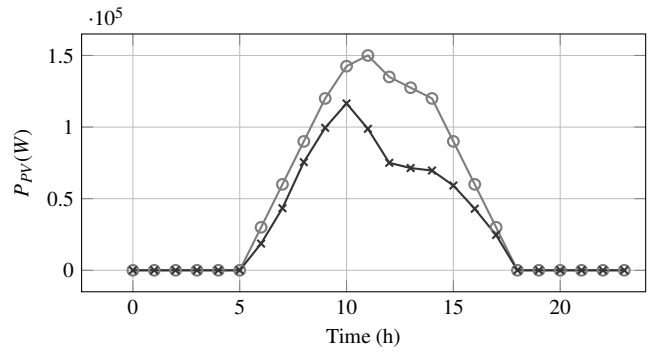
\begin{figure}[bth]
    \centering
    \footnotesize
    \input{Figures/SvSmax.tex}
    \caption{Optimal power generation (\ref{tikz:pv2}) and available power solar generation (\ref{tikz:pvmax2}) for Case 2.  The solar panels reduce their generation during the peak in order to maintain islanded operation.}
    \label{fig:SvSmax}
\end{figure}

\begin{figure}[bth]
    \centering
    \footnotesize
    \input{Figures/WvWmax.tex}
    \caption{Optimal wind generation (\ref{tikz:pw2}) and available wind generation (\ref{tikz:pwmax2}) for Case 2. The wind turbines reduce their generation maintaining the islanded operation an the balance between generation and power demand.}
    \label{fig:WvWmax}
\end{figure}
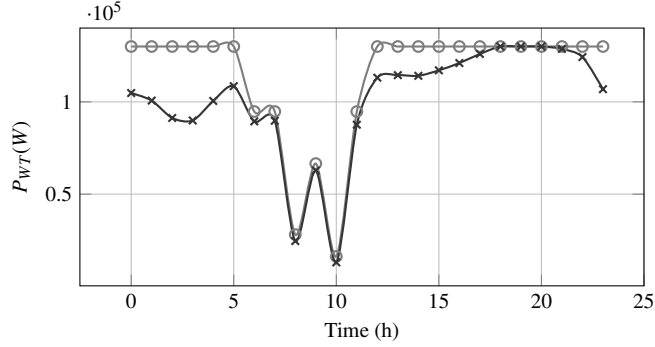

The photovoltaic generation uses almost all of its available resources, limiting its generation in periods of maximum irradiance, as depicted in Fig \ref{fig:SvSmax}. In this way, the microgrid keeps an equilibrium between generation and demand. Likewise, wind turbines reduce their generation to balance the load and generation. The behavior related to the wind turbines is shown in Fig \ref{fig:WvWmax}.

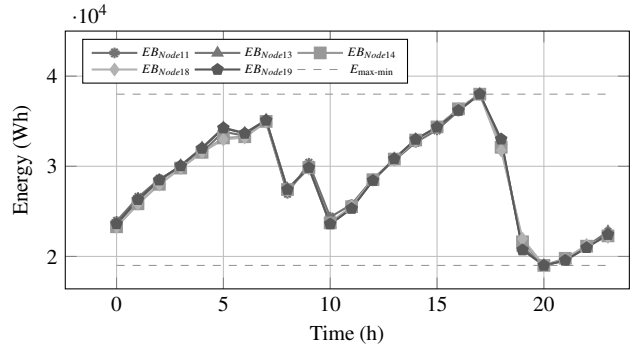
\begin{figure}[tbh]
    \centering
    \footnotesize
    \input{Figures/Baterias.tex}
    \caption{Energy stored in each battery for Case 2, considering the SOC of the batteries between $30\%$ and $100\%$ of the batteries energy.}
    \label{fig:Baterias}
\end{figure}

Eventually, the SOC and the batteries' initial state is the same as in the previous case. The batteries tend to charge up when existing a surplus of power supplied by the generation, seizing all the resources. However, the batteries try to inject power to the microgrid to maintain the balance in periods that the load power overtakes the generation, as shown in Fig \ref{fig:Baterias}; in these periods, the tertiary control requires to minimize the purchase of energy to the main grid and reduce costs in the peak of demand. 

\subsection{Software implementation}

All numerical experiments were performed in \textit{CvxPy}, a domain-specific language for convex optimization embedded in Python \cite{diamond2016cvxpy}. This module allows to write the model in a natural optimization syntax and call free and/or commercial packages to solve the problem. Different solvers were used to comparing the model performance in terms of operative cost and execution time.   Results are presented in Table \ref{tab:comparacion_solvers}, using a personal computer with 64-bit Operating system and processor Intel (R) Core (TM) i7-8700 CPU, 3.20 GHz, and 8GB of RAM. 

\begin{table*}
    \centering
    \caption{Performance comparison between different solvers.}
    \label{tab:comparacion_solvers}
    \begin{tabular}{cccccc}
    \toprule
        Case & Solver & Execution time (s) & Iterations & Operative cost ($\$$) & Status \\
    \midrule
          & Mosek   &    0.13  &    17   & -146340.52    &  Optimal  \\
     Case & ECOS    &	 0.19  &	33	 & -146340.52    &	Optimal  \\
     I    & ECOS BB	&    0.35  &	33	 & -146340.52    &	Optimal  \\
          & SCS	    &    1.55  &	2560 &	         	 &  Solved/Inaccurate \\
    \hline
          & Mosek   &    0.32  &    36   &  7887.18      & Optimal  \\
     Case & ECOS	&    0.12  &	21   &	7887.18      &	Optimal \\
     II   & ECOS BB	&    0.18  &	21	 &  7887.18	     &  Optimal \\
          & SCS	    &    3.1   &	121  &            	 & Solved/Inaccurate \\
    \bottomrule
    \end{tabular}
\end{table*}

Except for SCS, all solvers achieve the same optimal solution, although with different time calculations and number of iterations.  The model achieves the same global solution as was expected for convex optimization models.  Time calculation was less than 1s for the solvers in which the model achieves convergence.  This performance is important because a tertiary control must be executed in real-time, with the architecture already depicted in Fig \ref{fig:power_electronics}.

In practice, the optimization model is executed by a central controller that may be a personal computer or even a small single-board computer. Results presented in Table \ref{tab:comparacion_solvers} demonstrate the proposed model is suitable for the former case.  However, it is crucial to evaluate the model performance for a single-board computer.   

In this work, a Raspberry-pi is proposed as a central controller; this is a low-cost small single-board computer used in industrial applications. The computation power of such a small device is much lower compared to a personal computer.  However, it has some practical advantages, namely: the size of this board is 65mm x 30mm, so it can be integrated, for example, into the measurement unit of the microgrid; its cost is more than 20 times lower than a personal or industrial computer; it does not have, nor does it require, additional elements such as graphic or sound cards; in general, it can be used as an embedded application for central control.

All simulations presented in Section 6 were executed again in a Rasberry-pi model BV 1.2 with a 1.2 GHz 64-bit quad-core processor, onboard 802.11n WiFi, Bluetooth, and USB boot capabilities. This board uses a Linux operative system and includes a Python interpreter. This Rasberry-pi allowed the installation of the module \textit{CvxPy} for solving the optimization problems.  The solver ECOS was used since this solver was explicitly designed for embedded applications \cite{6669541}.   For Case 1, the same results as Section 6 were obtained 24.6 seconds; for Case 2, the same results were obtained in 33.93 seconds. In both cases, the response of the convex-TC implementation was suitable for real-time requirements mentioned in section \ref{sec:tertiary}.

\subsection{Extension to power distribution networks}
\begin{figure*}
    \centering
    \input{Figures/IEEE123}
    \caption{Modified IEEE123 test distribution system }
\label{fig:ieee123}
\end{figure*}
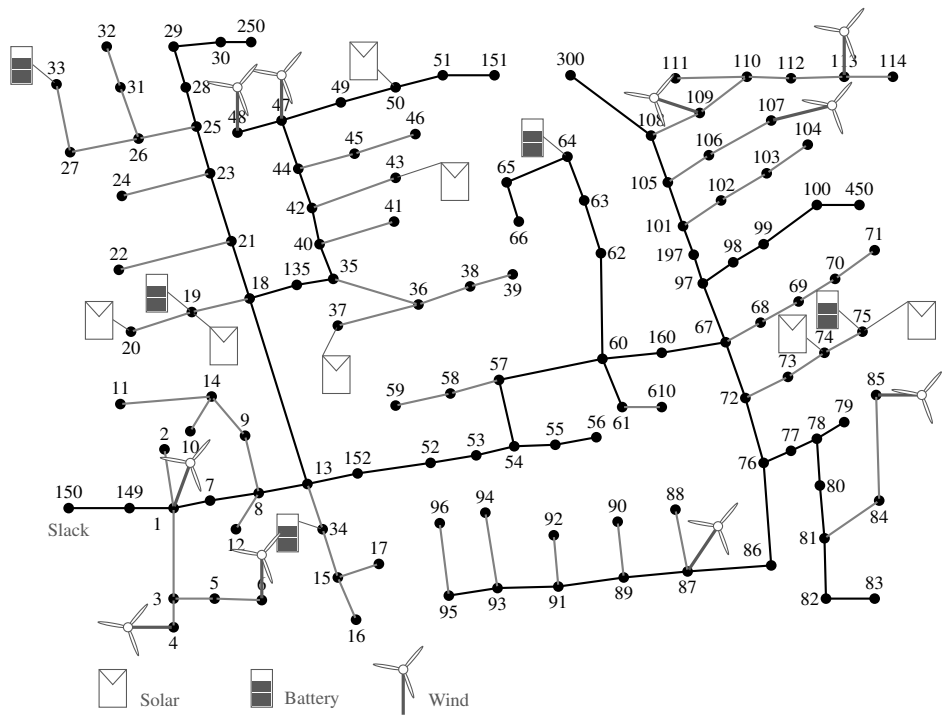

Although the proposed model was explicitly designed for microgrids, it may be extended to large power distribution systems; therefore, numerical experiments were performed in a modified version of the IEEE 123 node test feeder \cite{8063903}. This model includes overhead and underground single-phase lines, capacitor banks, unbalanced loads with constant current, constant impedance, and constant power models.  This test system was modified to include wind and solar energy and battery energy storage, as shown in Fig \ref{fig:ieee123}.

\begin{figure}[bth]
    \centering
    \footnotesize
    \input{Figures/Results2_IEEE123}
    \caption{General results for the IEEE 123 node test system.}
    \label{fig:resultados_ieee123}
\end{figure}
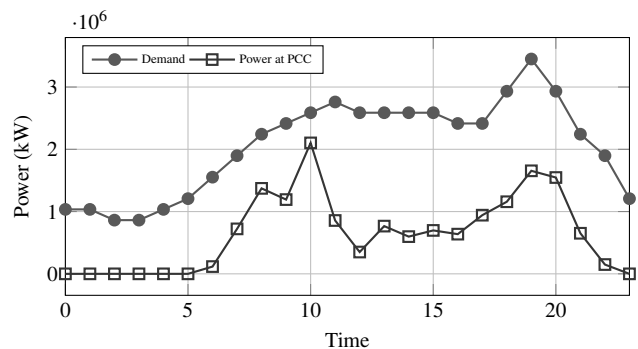

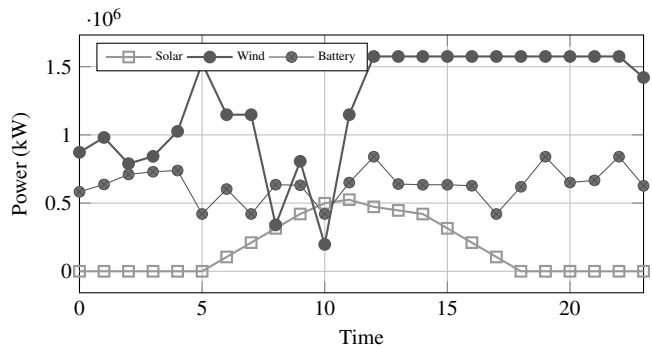
\begin{figure}[bth]
    \centering
    \footnotesize
    \input{Figures/Results_IEEE123}
    \caption{Total power generation for each distributed resource in the IEEE 123 node test system.}
    \label{fig:resultados_x_recurso_ieee123}
\end{figure}

The proposed algorithm was executed under the same considerations as Case 1, using a personal computer with a 64-bit Operating system and processor Intel (R) Core (TM) i7-8700 CPU, 3.20 GHz, and 8GB of RAM.  The model, in this case, consists of 197040 constraints, with 32736 convex cones and 196656 scalar variables. Despite its large size, the problem was solved in 2.26 seconds on average.  Figure \ref{fig:resultados_ieee123} shows the cycle of charge and discharge of battery energy storage devices.   Results for each distributed resource are depicted in Fig \ref{fig:resultados_x_recurso_ieee123}.

\section{Conclusions}

In this paper, an convex optimization model for the tertiary control of unbalanced microgrids was proposed. This model considered 24h operation with unbalanced loads; Also, it considered a realistic model of the loads, renewable energy sources, and storage devices. A Wirtinger linearization was used for the power flow equations, the exponential model of the loads, and the quadratic model of the battery energy storage. The method was implemented under changing wind and solar energy penetration. The proposed methodology ensured the optimum operation of the microgrid components. 

A real-time, low-cost, accurate implementation of the model was also considered. 
A Raspberry-pi was proposed as a central controller. The computation power of such a small device is much lower compared to a personal computer.  However, the proposed model was executed successfully in this device.  In this way, it was demonstrated the proposed model might be used for real-time operation. This implementation was a proof of concept under the worst-case scenario.

An intelligent predictor system based on prevision weather data (wind speed, solar irradiance, and temperature) has to be designed in further work. Prediction uncertainties on the renewable energies and load demand predictions have to be considered to design an autonomous system for an isolated area.

\section*{References}

\bibliographystyle{elsarticle-num}
\bibliography{Bibiografia}

\section*{Founding}
This research result is funded by UTP and the Colombian Science Ministry (Minciencias), project 111077657914, contract 031-2018, and project 321-2019.

\section*{Credit author statement}
\textit{Diego-Alejandro Ramirez}: methodology, software, writing. \textit{Alejandro Garc\'es}: conceptualization, methodology,  writing, validation. \textit{Juan-Jos\'e Mora-Fl\'orez}: validation, writing, editing.

\section*{Declaration of Competing Interest}
The authors declare that they have no known competing financial interests or personal relationships that could have appeared to influence the work reported in this paper.
\end{document}

%% file: Figures/Tikz_Convertidores.tex
\begin{tikzpicture}[x=1mm,y=1mm, thick]
% comunicaciones
\node at (77.5,0) [green!70!blue, draw, text width=40] (C) {central controller};

\node at (50,0) [green!70!blue, draw, text width=40] (F) {Forecast module};

\draw[-latex, green!70!blue] (F) -- (C);

\draw[dashed, green!70!blue] (5,20) |- (77.5,-35) -- (C)  -- (77.5,20);
\node[rotate=90] at (95,25) {Main grid};
\node[green!70!blue] at (50,-30) {communications};
% red
\draw[very thick] (25,-20) -- +(0,60);
\draw (25,25) -- (70,25);
\draw (85,25) -- (90,25);
\draw[very thick] (65,20) -- +(0,10);
\draw[very thick] (90,20) -- +(0,10);
\draw (75,25) circle (5);
\draw (80,25) circle (5);
% solar
\draw[fill=white] (0,0) rectangle +(10,10);
\draw (0,0) -- +(10,10);
\node at (3,7) {dc};
\node at (7,3) {ac};
\draw (0,5) -- +(-5,0);
\draw[fill=white] (-15,0) rectangle +(10,10);
\draw (-15,0) -- +(10,10);
\node at (-12,7) {dc};
\node at (-8,3) {dc};
\draw (-15,5) -- +(-7,0);
\draw (10,5) -- +(15,0);
\begin{scope}[ very thick, scale=0.7, xshift=-150, yshift=-5]
	\fill[green!40!blue] (0,0) -- +(5,21) -- +(26,21) -- +(21,0) -- cycle;
	\foreach \x in {0,3,...,21}  \draw[-,white] (\x,0) -- +(5,21);
	\foreach \y in {0,3,...,21}  \draw[-,white] (0.25*\y,\y) -- +(21,0);		
\end{scope}
\draw[-latex,blue!60!green] (-15,13) -- +(10,0) node[right] {$P_s$};
\draw[-latex,blue!60!green] (2,13) -- +(10,0) node[right] {$P_s,Q_s$};

% Turbina
\draw[fill=white] (0,20) rectangle +(10,10);
\draw (0,20) -- +(10,10);
\node at (3,27) {dc};
\node at (7,23) {ac};
\draw (0,25) -- +(-5,0);
\draw[fill=white] (-15,20) rectangle +(10,10);
\draw (-15,20) -- +(10,10);
\node at (-12,27) {ac};
\node at (-8,23) {dc};
\draw (-15,25) -- +(-15,0);
\draw (10,25) -- +(15,0);
\node(b) at (-30,40) {};
\draw[rotate = -15, fill=gray!30](b)+(5,0) ellipse (5 and 0.5);
\draw[rotate = 105, fill=gray!30](b)+(5,0) ellipse (5 and 0.5);
\draw[rotate = 225, fill=gray!30](b)+(5,0) ellipse (5 and 0.5);
\draw[-, fill=gray!30] (b)+(-0.5,0) -- +(-0.3,0) -- +(-0.6,-22) -- +(0.6,-22) -- +(0.3,0);
\draw[-, fill=gray!30] (b) circle (1);
\draw[-latex,blue!60!green] (-15,33) -- +(10,0) node[right] {$P_w$};
\draw[-latex,blue!60!green] (2,33) -- +(10,0) node[right] {$P_w,Q_w$};

% bateria
\draw[fill=white] (0,-20) rectangle +(10,10);
\draw (0,-20) -- +(10,10);
\node at (3,-13) {dc};
\node at (7,-17) {ac};
\draw (0,-15) -- +(-5,0);
\draw[fill=white] (-15,-20) rectangle +(10,10);
\draw (-15,-20) -- +(10,10);
\node at (-12,-13) {dc};
\node at (-8,-17) {dc};
\draw (-15,-15) -- +(-15,0);
\draw (10,-15) -- +(15,0);
\draw[-latex,blue!60!green] (-15,-7) -- +(10,0) node[right] {$P_B$};
\draw[-latex,blue!60!green] (0,-23) -- +(10,0) node[right] {$P_B,Q_B$};
\begin{scope}[xshift = -100, yshift = -55, scale=0.8]
	\fill [black!80] (3,0) rectangle +(14,15);
	\fill [orange] (3,15) rectangle +(14,-5);
	\fill [orange] (10,15) ellipse (7 and 2);
	\fill [orange] (10,10) ellipse (7 and 2);
	\fill [gray!40] (10,15) ellipse (6 and 1.5);
	\fill [black!80] (10,0) ellipse (7 and 2);
	\fill [gray] (7,15) rectangle +(6,2);
	\draw [gray, fill=gray!50, thick] (10,17) ellipse (3 and 1);
	\fill [gray] (10,15) ellipse (3 and 1);	
	\fill[white, opacity=0.2] (10,-2) -- +(-7,0) -- +(-7,17) -- cycle;
\end{scope}
\end{tikzpicture}

%% file: Figures/tikz_hypergraph.tex
\begin{tikzpicture}[x=1mm,y=1mm]
\draw[thick, gray, fill=blue!50!green!30] (0,10) ellipse (5 and 20);
\node at (0,-12) {hypernode k};

\draw[thick, gray, fill=blue!50!green!30] (50,10) ellipse (5 and 20);
\node at (50,-12) {hypernode m};

\draw[thick, gray, fill=blue!80!green!30] (27,12) ellipse (10 and 15);
\node at (27,-6) {hyperbranch km};

\draw[very thick] (0,0) -- +(0,5) node[right] {phase C};
\draw[very thick] (0,10) -- +(0,5) node[right] {phase B};
\draw[very thick] (0,20) -- +(0,5) node[right] {phase A};

\draw[very thick] (50,0) -- +(0,5) node[right] {phase C};
\draw[very thick] (50,10) -- +(0,5) node[right] {phase B};
\draw[very thick] (50,20) -- +(0,5) node[right] {phase A};

\draw[-latex](0,2.5) -- +(25,0);
\draw[-latex](0,12.5) -- +(25,0);
\draw[-latex](0,22.5) -- +(25,0);
\draw(25,2.5) -- +(25,0);
\draw(25,12.5) -- +(25,0);
\draw(25,22.5) -- +(25,0);
\end{tikzpicture}

%% file: Figures/tikz_cigre.tex
\centering
\footnotesize
\begin{tikzpicture}[x=1mm,y=1mm]
        \draw[very thick] (-10,100) node[left] {1} -- +(20,0) node[right] {400 V};
        \draw[thick] (0,100) -- (0,10);
        \draw[thick] (0,10) -- (-23,10);
        \draw[thick] (-23,80) -- (0,80);
        \draw[thick] (-27,70) -| (30,60);
        \draw[thick] (0,50) -- (-23,50);
        \draw[thick] (20,20) -- (0,20);
        \fill (0,90) circle (1) node[right, xshift=2] {2};
        \fill (0,80) circle (1) node[right, xshift=2] {3};
        \fill (0,70) circle (1) node[right, xshift=2, yshift=5] {4};
        \fill (0,60) circle (1) node[right, xshift=3] {5};
        \fill (0,50) circle (1) node[right, xshift=2, yshift=5] {6};
        \fill (0,40) circle (1) node[right, xshift=2] {7};
        \fill (0,30) circle (1) node[right, xshift=2] {8};
        \fill (0,20) circle (1) node[right, xshift=2, yshift = 5] {9};
        \fill (0,10) circle (1) node[right, xshift=2] {10};
        \fill (-23,80) circle (1) node[left, xshift=-2] {11};
        \fill (-23,70) circle (1) node[above, yshift=2] {12};
        \fill (10,70) circle (1) node[above, yshift=2] {15};
        \fill (20,70) circle (1) node[above, yshift=2] {16};
        \fill (30,70) circle (1) node[above, yshift=2] {17};
        \fill (30,60) circle (1) node[right, xshift=2] {18};
        \fill (-23,10) circle (1) node[below, yshift=-3] {14};
        \fill (-23,50) circle (1) node[below, xshift = 5, yshift=-2] {13};
        \fill (20,20) circle (1) node[above, yshift=3] {19};        
\node at (25,90) [text width = 100] {\scriptsize{Overhead line
$4\times 120 mm^2$ Al XLPE twisted cable, pole to pole distance $35m$}};
\node at (-11,82) {\tiny{$4\times 6mm^2$Cu, 30m}};
\node at (-11,72) {\tiny{$4\times 16mm^2$Cu, 30m}};
\node at (-12,52) {\tiny{$4\times 25mm^2$Cu, 30m}};
\node at (11,18) {\tiny{$4\times 6mm^2$Cu, 30m}};
\node at (-12,13) {\tiny{$4\times 16mm^2$Cu, 30m}};
\node at (15,68) {\tiny{$3\times 70mm^2$Al XLP, 56.4m}};

% usuario en nodo 11
\draw (-23,80) -- +(0,6);
\draw (-20,91) -- +(12,0);
\node at (-25,96) [text width=60]{\scriptsize{single residential consumer}};
\begin{scope}[xshift=-65, yshift=253, scale=0.2]
\draw[black,  fill=gray] (-22,2) -- +(4,15) -- +(40,15) -- +(44,0) -- cycle;
\draw[black,  fill=gray] (-8,17) -- +(0,3) -- +(4,3) -- +(4,0) -- cycle;
\draw[gray!30,top color=gray!30, bottom color=white] (-20,1) rectangle +(40,-2);
\draw[gray!30,top color=gray!30, bottom color=white] (-20,-1) rectangle +(40,-2);
\draw[gray!30,top color=gray!30, bottom color=white] (-20,-3) rectangle +(40,-2);
\draw[gray!30,top color=gray!30, bottom color=white] (-20,-5) rectangle +(40,-2);
\draw[gray!30,top color=gray!30, bottom color=white] (-20,-7) rectangle +(40,-2);
\draw[gray!30,top color=gray!30, bottom color=white] (-20,-9) rectangle +(40,-2);
\draw[gray!30,top color=gray!30, bottom color=white] (-20,-11) rectangle +(40,-2);
\draw[gray!30,top color=gray!30, bottom color=white] (-20,-13) rectangle +(40,-2);
\draw[black, fill=white] (-17,-9) rectangle +(6,7);
\draw[black] (-14,-9) -- +(0,7);
\draw[black] (-17,-5.5) -- +(6,0);
\draw[black, fill=white] (-8,-9) rectangle +(6,7);
\draw[black] (-5,-9) -- +(0,7);
\draw[black] (-8,-5.5) -- +(6,0);
\draw[black,fill=white] (11,-9) rectangle +(6,7);
\draw[black] (14,-9) -- +(0,7);
\draw[black] (11,-5.5) -- +(6,0);
\draw[black, fill=gray] (1,-15) rectangle +(7,14);
\draw[gray!30,fill] (7,-8) circle (0.5);
\draw[black] (-20,1) rectangle +(40,-16);
\draw[black!70, fill] (-22,1) rectangle +(44,1);
\draw[black!70, fill] (0,-15.5) rectangle +(9,1);
\end{scope}
% energy storage device
\draw [thick] (-33,67) rectangle +(6,6);
\draw [thick] (-33,67) -- +(6,6);
\node at (-31.5,71.5) {\tiny{ac}};
\node at (-28.5,69) {\tiny{ac}};
% energy renewable resourse
\draw [thick] (-33,33) rectangle +(6,6);
\draw [thick] (-33,33) -- +(6,6);
\node at (-31.5,37.5) {\tiny{dc}};
\node at (-28.5,35) {\tiny{ac}};
\draw[thin] (-27,36) -| +(4,8);
% cargas del nodo 13
\draw (-40,44) rectangle +(17,10);
\node at (-30,59) {\scriptsize{Group of 4 residences}};
\begin{scope}[gray, xshift=-85, yshift=150, scale=0.2]
\draw[black,  fill=gray] (-22,2) -- +(4,15) -- +(40,15) -- +(44,0) -- cycle;
\draw[black,  fill=gray] (-8,17) -- +(0,3) -- +(4,3) -- +(4,0) -- cycle;
\draw[gray!30,top color=gray!30, bottom color=white] (-20,1) rectangle +(40,-2);
\draw[gray!30,top color=gray!30, bottom color=white] (-20,-1) rectangle +(40,-2);
\draw[gray!30,top color=gray!30, bottom color=white] (-20,-3) rectangle +(40,-2);
\draw[gray!30,top color=gray!30, bottom color=white] (-20,-5) rectangle +(40,-2);
\draw[gray!30,top color=gray!30, bottom color=white] (-20,-7) rectangle +(40,-2);
\draw[gray!30,top color=gray!30, bottom color=white] (-20,-9) rectangle +(40,-2);
\draw[gray!30,top color=gray!30, bottom color=white] (-20,-11) rectangle +(40,-2);
\draw[gray!30,top color=gray!30, bottom color=white] (-20,-13) rectangle +(40,-2);
\draw[black, fill=white] (-17,-9) rectangle +(6,7);
\draw[black] (-14,-9) -- +(0,7);
\draw[black] (-17,-5.5) -- +(6,0);
\draw[black, fill=white] (-8,-9) rectangle +(6,7);
\draw[black] (-5,-9) -- +(0,7);
\draw[black] (-8,-5.5) -- +(6,0);
\draw[black,fill=white] (11,-9) rectangle +(6,7);
\draw[black] (14,-9) -- +(0,7);
\draw[black] (11,-5.5) -- +(6,0);
\draw[black, fill=gray] (1,-15) rectangle +(7,14);
\draw[gray!30,fill] (7,-8) circle (0.5);
\draw[black] (-20,1) rectangle +(40,-16);
\draw[black!70, fill] (-22,1) rectangle +(44,1);
\draw[black!70, fill] (0,-15.5) rectangle +(9,1);
\end{scope}
\begin{scope}[gray, xshift=-85, yshift=128, scale=0.2]
\draw[black,  fill=gray] (-22,2) -- +(4,15) -- +(40,15) -- +(44,0) -- cycle;
\draw[black,  fill=gray] (-8,17) -- +(0,3) -- +(4,3) -- +(4,0) -- cycle;
\draw[gray!30,top color=gray!30, bottom color=white] (-20,1) rectangle +(40,-2);
\draw[gray!30,top color=gray!30, bottom color=white] (-20,-1) rectangle +(40,-2);
\draw[gray!30,top color=gray!30, bottom color=white] (-20,-3) rectangle +(40,-2);
\draw[gray!30,top color=gray!30, bottom color=white] (-20,-5) rectangle +(40,-2);
\draw[gray!30,top color=gray!30, bottom color=white] (-20,-7) rectangle +(40,-2);
\draw[gray!30,top color=gray!30, bottom color=white] (-20,-9) rectangle +(40,-2);
\draw[gray!30,top color=gray!30, bottom color=white] (-20,-11) rectangle +(40,-2);
\draw[gray!30,top color=gray!30, bottom color=white] (-20,-13) rectangle +(40,-2);
\draw[black, fill=white] (-17,-9) rectangle +(6,7);
\draw[black] (-14,-9) -- +(0,7);
\draw[black] (-17,-5.5) -- +(6,0);
\draw[black, fill=white] (-8,-9) rectangle +(6,7);
\draw[black] (-5,-9) -- +(0,7);
\draw[black] (-8,-5.5) -- +(6,0);
\draw[black,fill=white] (11,-9) rectangle +(6,7);
\draw[black] (14,-9) -- +(0,7);
\draw[black] (11,-5.5) -- +(6,0);
\draw[black, fill=gray] (1,-15) rectangle +(7,14);
\draw[gray!30,fill] (7,-8) circle (0.5);
\draw[black] (-20,1) rectangle +(40,-16);
\draw[black!70, fill] (-22,1) rectangle +(44,1);
\draw[black!70, fill] (0,-15.5) rectangle +(9,1);
\end{scope}
\begin{scope}[gray, xshift=-113, yshift=128, scale=0.2]
\draw[black,  fill=gray] (-22,2) -- +(4,15) -- +(40,15) -- +(44,0) -- cycle;
\draw[black,  fill=gray] (-8,17) -- +(0,3) -- +(4,3) -- +(4,0) -- cycle;
\draw[gray!30,top color=gray!30, bottom color=white] (-20,1) rectangle +(40,-2);
\draw[gray!30,top color=gray!30, bottom color=white] (-20,-1) rectangle +(40,-2);
\draw[gray!30,top color=gray!30, bottom color=white] (-20,-3) rectangle +(40,-2);
\draw[gray!30,top color=gray!30, bottom color=white] (-20,-5) rectangle +(40,-2);
\draw[gray!30,top color=gray!30, bottom color=white] (-20,-7) rectangle +(40,-2);
\draw[gray!30,top color=gray!30, bottom color=white] (-20,-9) rectangle +(40,-2);
\draw[gray!30,top color=gray!30, bottom color=white] (-20,-11) rectangle +(40,-2);
\draw[gray!30,top color=gray!30, bottom color=white] (-20,-13) rectangle +(40,-2);
\draw[black, fill=white] (-17,-9) rectangle +(6,7);
\draw[black] (-14,-9) -- +(0,7);
\draw[black] (-17,-5.5) -- +(6,0);
\draw[black, fill=white] (-8,-9) rectangle +(6,7);
\draw[black] (-5,-9) -- +(0,7);
\draw[black] (-8,-5.5) -- +(6,0);
\draw[black,fill=white] (11,-9) rectangle +(6,7);
\draw[black] (14,-9) -- +(0,7);
\draw[black] (11,-5.5) -- +(6,0);
\draw[black, fill=gray] (1,-15) rectangle +(7,14);
\draw[gray!30,fill] (7,-8) circle (0.5);
\draw[black] (-20,1) rectangle +(40,-16);
\draw[black!70, fill] (-22,1) rectangle +(44,1);
\draw[black!70, fill] (0,-15.5) rectangle +(9,1);
\end{scope}
\begin{scope}[gray, xshift=-113, yshift=150, scale=0.2]
\draw[black,  fill=gray] (-22,2) -- +(4,15) -- +(40,15) -- +(44,0) -- cycle;
\draw[black,  fill=gray] (-8,17) -- +(0,3) -- +(4,3) -- +(4,0) -- cycle;
\draw[gray!30,top color=gray!30, bottom color=white] (-20,1) rectangle +(40,-2);
\draw[gray!30,top color=gray!30, bottom color=white] (-20,-1) rectangle +(40,-2);
\draw[gray!30,top color=gray!30, bottom color=white] (-20,-3) rectangle +(40,-2);
\draw[gray!30,top color=gray!30, bottom color=white] (-20,-5) rectangle +(40,-2);
\draw[gray!30,top color=gray!30, bottom color=white] (-20,-7) rectangle +(40,-2);
\draw[gray!30,top color=gray!30, bottom color=white] (-20,-9) rectangle +(40,-2);
\draw[gray!30,top color=gray!30, bottom color=white] (-20,-11) rectangle +(40,-2);
\draw[gray!30,top color=gray!30, bottom color=white] (-20,-13) rectangle +(40,-2);
\draw[black, fill=white] (-17,-9) rectangle +(6,7);
\draw[black] (-14,-9) -- +(0,7);
\draw[black] (-17,-5.5) -- +(6,0);
\draw[black, fill=white] (-8,-9) rectangle +(6,7);
\draw[black] (-5,-9) -- +(0,7);
\draw[black] (-8,-5.5) -- +(6,0);
\draw[black,fill=white] (11,-9) rectangle +(6,7);
\draw[black] (14,-9) -- +(0,7);
\draw[black] (11,-5.5) -- +(6,0);
\draw[black, fill=gray] (1,-15) rectangle +(7,14);
\draw[gray!30,fill] (7,-8) circle (0.5);
\draw[black] (-20,1) rectangle +(40,-16);
\draw[black!70, fill] (-22,1) rectangle +(44,1);
\draw[black!70, fill] (0,-15.5) rectangle +(9,1);
\end{scope}
% cargas nodo 19
\draw [thick] (25,13) rectangle +(6,6);
\draw [thick] (25,13) -- +(6,6);
\node at (26.5,17.5) {\tiny{ac}};
\node at (29,14.5) {\tiny{dc}};
\draw[thin] (-27,36) -| +(4,8);
%\node at (23,11) {\scriptsize{photovoltaics+storage}};
\draw (20,20) -| +(3,-4) -- +(5,-4);
\draw (20,20) -| +(3,4) -- +(5,4);
\node at (25,31) [text width=70]{\scriptsize{single residential consumer}};
\begin{scope}[gray, xshift=80, yshift=75, scale=0.2]
\draw[black,  fill=gray] (-22,2) -- +(4,15) -- +(40,15) -- +(44,0) -- cycle;
\draw[black,  fill=gray] (-8,17) -- +(0,3) -- +(4,3) -- +(4,0) -- cycle;
\draw[gray!30,top color=gray!30, bottom color=white] (-20,1) rectangle +(40,-2);
\draw[gray!30,top color=gray!30, bottom color=white] (-20,-1) rectangle +(40,-2);
\draw[gray!30,top color=gray!30, bottom color=white] (-20,-3) rectangle +(40,-2);
\draw[gray!30,top color=gray!30, bottom color=white] (-20,-5) rectangle +(40,-2);
\draw[gray!30,top color=gray!30, bottom color=white] (-20,-7) rectangle +(40,-2);
\draw[gray!30,top color=gray!30, bottom color=white] (-20,-9) rectangle +(40,-2);
\draw[gray!30,top color=gray!30, bottom color=white] (-20,-11) rectangle +(40,-2);
\draw[gray!30,top color=gray!30, bottom color=white] (-20,-13) rectangle +(40,-2);
\draw[black, fill=white] (-17,-9) rectangle +(6,7);
\draw[black] (-14,-9) -- +(0,7);
\draw[black] (-17,-5.5) -- +(6,0);
\draw[black, fill=white] (-8,-9) rectangle +(6,7);
\draw[black] (-5,-9) -- +(0,7);
\draw[black] (-8,-5.5) -- +(6,0);
\draw[black,fill=white] (11,-9) rectangle +(6,7);
\draw[black] (14,-9) -- +(0,7);
\draw[black] (11,-5.5) -- +(6,0);
\draw[black, fill=gray] (1,-15) rectangle +(7,14);
\draw[gray!30,fill] (7,-8) circle (0.5);
\draw[black] (-20,1) rectangle +(40,-16);
\draw[black!70, fill] (-22,1) rectangle +(44,1);
\draw[black!70, fill] (0,-15.5) rectangle +(9,1);
\end{scope}
\draw (31,16) -- +(3,0);
% Bateria en el nodo 19
\begin{scope}[xshift = 100, yshift = 50, scale=0.4]
	\fill [black!80] (3,0) rectangle +(14,15);
	\fill [orange] (3,15) rectangle +(14,-5);
	\fill [orange] (10,15) ellipse (7 and 2);
	\fill [orange] (10,10) ellipse (7 and 2);
	\fill [gray!40] (10,15) ellipse (6 and 1.5);
	\fill [black!80] (10,0) ellipse (7 and 2);
	\fill [gray] (7,15) rectangle +(6,2);
	\draw [gray, fill=gray!50, thick] (10,17) ellipse (3 and 1);
	\fill [gray] (10,15) ellipse (3 and 1);	
	\fill[white, opacity=0.2] (10,-2) -- +(-7,0) -- +(-7,17) -- cycle;
\end{scope}
% Panel solar nodo 19
\begin{scope}[ scale=0.3, xshift=310, yshift=118]	
	\fill[gray!50] (0,0) -- +(5,5) -- +(32,7) -- +(21,0) -- cycle;
	\fill[green!40!blue] (0,0) -- +(5,21) -- +(26,21) -- +(21,0) -- cycle;
	\foreach \x in {0,3,...,21}  \draw[-,white] (\x,0) -- +(5,21);
	\foreach \y in {0,3,...,21}  \draw[-,white] (0.25*\y,\y) -- +(21,0);
\end{scope}
% cargas en el nodo 18
\draw (26,54) -| +(4,5);
\draw [thick] (20,40) rectangle +(6,6);
\draw [thick] (20,40) -- +(6,6);
\node at (22,44) {\tiny{dc}};
\node at (24,41) {\tiny{ac}};
\draw (30,54) |- +(-4,-11);
\draw (20,43) -- +(-5,0);
% apartment building
\node at (12,57) [text width=30] {\scriptsize{apartment building}}; 
\begin{scope}[scale = 0.3, xshift=180, yshift=470]    
  \draw[fill=gray] (0,0) rectangle +(28,16);    
  \draw[fill=gray!50] (5,50) rectangle +(5,4);
  \draw[fill=white] (-0.5,48) rectangle +(29,2); % corniza
  % tercer piso
  \draw[fill=gray] (0,32) rectangle +(28,16);
  \draw[black, fill=white] (2,37) rectangle +(6,8);
  \draw[black] (5,37) -- +(0,8); 
  \draw[black] (2,41) -- +(6,0); 
  \draw[black, fill=white] (1.5,36) rectangle +(7,1);
  
  \draw[black, fill=white] (11,37) rectangle +(6,8);
  \draw[black] (14,37) -- +(0,8); 
  \draw[black] (11,41) -- +(6,0); 
  \draw[black, fill=white] (10.5,36) rectangle +(7,1);
  
  \draw[black, fill=white] (20,37) rectangle +(6,8);
  \draw[black] (23,37) -- +(0,8); 
  \draw[black] (20,41) -- +(6,0); 
  \draw[black, fill=white] (19.5,36) rectangle +(7,1);
  % segundo piso
  \draw[fill=gray] (0,16) rectangle +(28,16);
  \draw[black, fill=white] (2,21) rectangle +(6,8);
  \draw[black] (5,21) -- +(0,8); 
  \draw[black] (2,25) -- +(6,0); 
  \draw[black, fill=white] (1.5,20) rectangle +(7,1);
  
  \draw[black, fill=white] (11,21) rectangle +(6,8);
  \draw[black] (14,21) -- +(0,8); 
  \draw[black] (11,25) -- +(6,0); 
  \draw[black, fill=white] (10.5,20) rectangle +(7,1);
  
  \draw[black, fill=white] (20,21) rectangle +(6,8);
  \draw[black] (23,21) -- +(0,8); 
  \draw[black] (20,25) -- +(6,0); 
  \draw[black, fill=white] (19.5,20) rectangle +(7,1);
  % primer piso
  \draw[black, fill=white] (2,5) rectangle +(6,8);
  \draw[black] (5,5) -- +(0,8); 
  \draw[black] (2,9) -- +(6,0); 
  \draw[black, fill=white] (1.5,4) rectangle +(7,1);
  
  \draw[black, fill=white] (11,5) rectangle +(6,8);
  \draw[black] (14,5) -- +(0,8); 
  \draw[black] (11,9) -- +(6,0); 
  \draw[black, fill=white] (10.5,4) rectangle +(7,1);
 
  % puerta
     \draw[black, fill=white] (20,12) arc (180:0:2.5);
     \draw[black, fill=white] (20,5) rectangle +(5,8);  
     \fill (21,9) circle (0.5);
     \draw[black, fill=white] (19,3) rectangle +(7,2);
     \draw[black, fill=white] (19,1) rectangle +(7,2);
     \draw[black, fill=white] (19,-1) rectangle +(7,2);
     \draw[black, fill=white] (18,-1) rectangle +(2,7);
     \draw[black, fill=white] (25,-1) rectangle +(2,7);
\end{scope}
% cargas en el nodo 14
% apartment building
\node at (-27,4) {apartment building};
\begin{scope}[scale = 0.3, xshift=-410, yshift=70]    
  \draw[fill=gray] (0,0) rectangle +(28,16);    
  \draw[fill=gray!50] (5,50) rectangle +(5,4);
  \draw[fill=white] (-0.5,48) rectangle +(29,2); % corniza
  
  % tercer piso
  \draw[fill=gray] (0,32) rectangle +(28,16);
  \draw[black, fill=white] (2,37) rectangle +(6,8);
  \draw[black] (5,37) -- +(0,8); 
  \draw[black] (2,41) -- +(6,0); 
  \draw[black, fill=white] (1.5,36) rectangle +(7,1);
  
  \draw[black, fill=white] (11,37) rectangle +(6,8);
  \draw[black] (14,37) -- +(0,8); 
  \draw[black] (11,41) -- +(6,0); 
  \draw[black, fill=white] (10.5,36) rectangle +(7,1);
  
  \draw[black, fill=white] (20,37) rectangle +(6,8);
  \draw[black] (23,37) -- +(0,8); 
  \draw[black] (20,41) -- +(6,0); 
  \draw[black, fill=white] (19.5,36) rectangle +(7,1);
  % segundo piso
  \draw[fill=gray] (0,16) rectangle +(28,16);
  \draw[black, fill=white] (2,21) rectangle +(6,8);
  \draw[black] (5,21) -- +(0,8); 
  \draw[black] (2,25) -- +(6,0); 
  \draw[black, fill=white] (1.5,20) rectangle +(7,1);
  
  \draw[black, fill=white] (11,21) rectangle +(6,8);
  \draw[black] (14,21) -- +(0,8); 
  \draw[black] (11,25) -- +(6,0); 
  \draw[black, fill=white] (10.5,20) rectangle +(7,1);
  
  \draw[black, fill=white] (20,21) rectangle +(6,8);
  \draw[black] (23,21) -- +(0,8); 
  \draw[black] (20,25) -- +(6,0); 
  \draw[black, fill=white] (19.5,20) rectangle +(7,1);
  % primer piso
  \draw[black, fill=white] (2,5) rectangle +(6,8);
  \draw[black] (5,5) -- +(0,8); 
  \draw[black] (2,9) -- +(6,0); 
  \draw[black, fill=white] (1.5,4) rectangle +(7,1);
  
  \draw[black, fill=white] (11,5) rectangle +(6,8);
  \draw[black] (14,5) -- +(0,8); 
  \draw[black] (11,9) -- +(6,0); 
  \draw[black, fill=white] (10.5,4) rectangle +(7,1);
 
  % puerta
     \draw[black, fill=white] (20,12) arc (180:0:2.5);
     \draw[black, fill=white] (20,5) rectangle +(5,8);  
     \fill (21,9) circle (0.5);
     \draw[black, fill=white] (19,3) rectangle +(7,2);
     \draw[black, fill=white] (19,1) rectangle +(7,2);
     \draw[black, fill=white] (19,-1) rectangle +(7,2);
     \draw[black, fill=white] (18,-1) rectangle +(2,7);
     \draw[black, fill=white] (25,-1) rectangle +(2,7);
\end{scope}

% Turbina Nodo 12
\draw (-33,70) -- +(-5,0);
\node(b) at (-38,88) {};
\draw[rotate = -15, fill=gray!30](b)+(5,0) ellipse (5 and 0.5);
\draw[rotate = 105, fill=gray!30](b)+(5,0) ellipse (5 and 0.5);
\draw[rotate = 225, fill=gray!30](b)+(5,0) ellipse (5 and 0.5);
\draw[-, fill=gray!30] (b)+(-0.5,0) -- +(-0.3,0) -- +(-0.6,-22) -- +(0.6,-22) -- +(0.3,0);
\draw[-, fill=gray!30] (b) circle (1);

% Turbina Nodo 13
\draw (-23,38) |- (-12,23);
\draw [thick, fill=white] (-20,20) rectangle +(6,6);
\draw [thick] (-20,20) -- +(6,6);
\node at (-18,24) {\tiny{ac}};
\node at (-16,22) {\tiny{ac}};

\node(b) at (-12,38) {};
\draw[rotate = -15, fill=gray!30](b)+(5,0) ellipse (5 and 0.5);
\draw[rotate = 105, fill=gray!30](b)+(5,0) ellipse (5 and 0.5);
\draw[rotate = 225, fill=gray!30](b)+(5,0) ellipse (5 and 0.5);
\draw[-, fill=gray!30] (b)+(-0.5,0) -- +(-0.3,0) -- +(-0.6,-22) -- +(0.6,-22) -- +(0.3,0);
\draw[-, fill=gray!30] (b) circle (1);

% Bateria Nodo 13
\draw (-33,36) -- +(-10,0);
\begin{scope}[xshift = -128, yshift = 90, scale=0.4]
	\fill [black!80] (3,0) rectangle +(14,15);
	\fill [orange] (3,15) rectangle +(14,-5);
	\fill [orange] (10,15) ellipse (7 and 2);
	\fill [orange] (10,10) ellipse (7 and 2);
	\fill [gray!40] (10,15) ellipse (6 and 1.5);
	\fill [black!80] (10,0) ellipse (7 and 2);
	\fill [gray] (7,15) rectangle +(6,2);
	\draw [gray, fill=gray!50, thick] (10,17) ellipse (3 and 1);
	\fill [gray] (10,15) ellipse (3 and 1);	
	\fill[white, opacity=0.2] (10,-2) -- +(-7,0) -- +(-7,17) -- cycle;
\end{scope}

% Panel nodo 13
\draw (-35,36) |- +(-2,-5);
\begin{scope}[scale=0.3, xshift=-420, yshift=255]	
	\fill[gray!50] (0,0) -- +(5,5) -- +(32,7) -- +(21,0) -- cycle;
	\fill[green!40!blue] (0,0) -- +(5,21) -- +(26,21) -- +(21,0) -- cycle;
	\foreach \x in {0,3,...,21}  \draw[-,white] (\x,0) -- +(5,21);
	\foreach \y in {0,3,...,21}  \draw[-,white] (0.25*\y,\y) -- +(21,0);
\end{scope}

% Bateria nodo 18
\begin{scope}[xshift = 18, yshift = 122, scale=0.4]
	\fill [black!80] (3,0) rectangle +(14,15);
	\fill [orange] (3,15) rectangle +(14,-5);
	\fill [orange] (10,15) ellipse (7 and 2);
	\fill [orange] (10,10) ellipse (7 and 2);
	\fill [gray!40] (10,15) ellipse (6 and 1.5);
	\fill [black!80] (10,0) ellipse (7 and 2);
	\fill [gray] (7,15) rectangle +(6,2);
	\draw [gray, fill=gray!50, thick] (10,17) ellipse (3 and 1);
	\fill [gray] (10,15) ellipse (3 and 1);	
	\fill[white, opacity=0.2] (10,-2) -- +(-7,0) -- +(-7,17) -- cycle;
\end{scope}

% Panel nodo 18
\begin{scope}[scale=0.3, xshift=90, yshift=370]	
	\fill[gray!50] (0,0) -- +(5,5) -- +(32,7) -- +(21,0) -- cycle;
	\fill[green!40!blue] (0,0) -- +(5,21) -- +(26,21) -- +(21,0) -- cycle;
	\foreach \x in {0,3,...,21}  \draw[-,white] (\x,0) -- +(5,21);
	\foreach \y in {0,3,...,21}  \draw[-,white] (0.25*\y,\y) -- +(21,0);
\end{scope}

% Bateria nodo 11
\begin{scope}[xshift = -28, yshift = 250, scale=0.4]
	\fill [black!80] (3,0) rectangle +(14,15);
	\fill [orange] (3,15) rectangle +(14,-5);
	\fill [orange] (10,15) ellipse (7 and 2);
	\fill [orange] (10,10) ellipse (7 and 2);
	\fill [gray!40] (10,15) ellipse (6 and 1.5);
	\fill [black!80] (10,0) ellipse (7 and 2);
	\fill [gray] (7,15) rectangle +(6,2);
	\draw [gray, fill=gray!50, thick] (10,17) ellipse (3 and 1);
	\fill [gray] (10,15) ellipse (3 and 1);	
	\fill[white, opacity=0.2] (10,-2) -- +(-7,0) -- +(-7,17) -- cycle;
\end{scope}
\draw [thick,fill=white] (-16,88) rectangle +(6,6);
\draw [thick] (-16,88) -- +(6,6);
\node at (-14,92.5) {\tiny{ac}};
\node at (-12,89.5) {\tiny{dc}};

% Panel nodo 12
\begin{scope}[scale=0.3, xshift=-100, yshift=580]	
	\fill[gray!50] (0,0) -- +(5,5) -- +(32,7) -- +(21,0) -- cycle;
	\fill[green!40!blue] (0,0) -- +(5,21) -- +(26,21) -- +(21,0) -- cycle;
	\foreach \x in {0,3,...,21}  \draw[-,white] (\x,0) -- +(5,21);
	\foreach \y in {0,3,...,21}  \draw[-,white] (0.25*\y,\y) -- +(21,0);
\end{scope}
\draw (-23,70) |- +(15,-5);
\draw [thick,fill=white] (-20,62) rectangle +(6,6);
\draw [thick] (-20,62) -- +(6,6);
\node at (-18.5,65.5) {\tiny{ac}};
\node at (-16,63.5) {\tiny{dc}};

% Bateria nodo 14
\begin{scope}[xshift = -100, yshift = 68, scale=0.4]
	\fill [black!80] (3,0) rectangle +(14,15);
	\fill [orange] (3,15) rectangle +(14,-5);
	\fill [orange] (10,15) ellipse (7 and 2);
	\fill [orange] (10,10) ellipse (7 and 2);
	\fill [gray!40] (10,15) ellipse (6 and 1.5);
	\fill [black!80] (10,0) ellipse (7 and 2);
	\fill [gray] (7,15) rectangle +(6,2);
	\draw [gray, fill=gray!50, thick] (10,17) ellipse (3 and 1);
	\fill [gray] (10,15) ellipse (3 and 1);	
	\fill[white, opacity=0.2] (10,-2) -- +(-7,0) -- +(-7,17) -- cycle;
\end{scope}
% Panel nodo 14
\begin{scope}[scale=0.3, xshift=-310, yshift=190]	
	\fill[gray!50] (0,0) -- +(5,5) -- +(32,7) -- +(21,0) -- cycle;
	\fill[green!40!blue] (0,0) -- +(5,21) -- +(26,21) -- +(21,0) -- cycle;
	\foreach \x in {0,3,...,21}  \draw[-,white] (\x,0) -- +(5,21);
	\foreach \y in {0,3,...,21}  \draw[-,white] (0.25*\y,\y) -- +(21,0);
\end{scope}
\draw (-29,10) -- +(0,10);
\draw [thick, fill=white] (-32,12) rectangle +(6,6);
\draw [thick] (-32,12) -- +(6,6);
\node at (-30,16.5) {\tiny{dc}};
\node at (-27.5,14) {\tiny{ac}};
\draw (-20,10) -- (-35,10);

\end{tikzpicture}

%% file: Figures/Costo.tex
\begin{tikzpicture}
\begin{axis}[width=9cm, height=5cm, xlabel={Time (h)},ylabel={Cost (\$/kW)},xmin=1,xmax=24,xmajorgrids,ymajorgrids]
\addplot [mark=o,color=blue!80!green, thick]
coordinates{
    (0,0.10165)(1,0.10165)(2,0.10165)(3,0.10115)(4,0.10165)(5,0.10165)(6,0.10165)(7,0.10165)(8,0.11322)(9,0.11322)(10,0.11322)(11,0.13236)(12,0.13236)(13,0.11322)(14,0.11322)(15,0.11322)(16,0.11322)(17,0.11322)(18,0.13236)(19,0.13236)(20,0.11322)(21,0.11322)(22,0.11322)(23,0.10165)
    };    
\label{fig:Cost}
\end{axis}
\end{tikzpicture}

%% file: Figures/GenvsDem_PowSale.tex
\begin{tikzpicture}
\begin{axis}[width=9cm, height=5cm, xlabel={Time (h)},ylabel={Power (kW)},xmin=1,xmax=24,xmajorgrids,ymajorgrids]
\addplot [mark=o,color=green!80!blue, thick]
coordinates{
    (0,56.070)(1,56.070)(2,46.725)(3,46.725)(4,56.070)(5,65.415)(6,84.105)(7,102.795)(8,121.485)(9,130.830)(10,140.175)(11,149.520)(12,140.175)(13,140.175)(14,140.175)(15,140.175)(16,130.830)(17,130.830)(18,158.865)(19,186.900)(20,158.865)(21,121.485)(22,102.795)(23,65.415)
    };    
\label{tikz:demand}

\addplot[mark=triangle, color=orange,  thick] coordinates {(0,122.9503)(1,122.9486)(2,126.4752)(3,116.5326)(4,122.9466)(5,117.9089)(6,102.5687)(7,122.4481)(8,117.5751)(9,186.0531)(10,156.7285)(11,287.7452)(12,312.8623)(13,235.1847)(14,227.6843)(15,197.6842)(16,172.6122)(17,142.6125)(18,182.6464)(19,168.17)(20,127.9723)(21,129.4938)(22,129.4947)(23,128.9885)};
\label{tikz:per}

\addplot[mark=square, color=blue!80, thick] coordinates {(0,-67.0843)(1,-67.0827)(2,-79.7861)(3,-69.9436)(4,-67.0807)(5,-52.956)(6,-19.7782)(7,-21.2148)(8,2.2201)(9,-55.4438)(10,-17.1594)(11,-135.1689)(12,-169.1228)(13,-94.0718)(14,-86.8847)(15,-58.0295)(16,-42.8153)(17,-13.7382)(18,-25.5296)(19,14.8837)(20,27.1046)(21,-9.9514)(22,-27.9493)(23,-63.8642)};
\label{tikz:slack}
\end{axis}
\end{tikzpicture}

%% file: Figures/SvSmax_PowSale.tex
\begin{tikzpicture}
\begin{axis}[width=9cm, height=5cm,xlabel={Time (h)},ylabel={$P_{PV}(W)$},xmax=24,legend pos=north east, xmajorgrids,ymajorgrids,legend style={nodes={scale=1, transform shape}}]
\addplot [smooth,color=green!80!blue,mark=o,thick]
coordinates{
    (0,0)(1,0)(2,0)(3,0)(4,0)(5,0)(6,30000)(7,60000)(8,90000)(9,120000)(10,142500)(11,150000)(12,135000)(13,127500)(14,120000)(15,90000)(16,60000)(17,30000)(18,0)(19,0)(20,0)(21,0)(22,0)(23,0)
    };    
%\addlegendentry{$P_{PV}max(t)$}
\label{tikz:pvmax}
\addplot[smooth,mark=x, color=blue!80!green, thick] coordinates {(0,0)(1,0)(2,0)(3,0)(4,0)(5,0)(6,29999.9926)(7,59999.9926)(8,89999.9933)(9,119999.9932)(10,142499.9931)(11,149991.9909)(12,134999.9941)(13,127499.9932)(14,119999.9932)(15,89999.9932)(16,59999.9933)(17,29999.9933)(18,0)(19,0)(20,0)(21,0)(22,0)(23,0)};
\label{tikz:ppv}
%\addlegendentry{$P_{PV}(t)$}

\end{axis}
\end{tikzpicture}

%% file: Figures/WvWmax_PowSale.tex
\begin{tikzpicture}
\begin{axis}[width=9cm, height=5.0cm, xlabel={time (h)},ylabel={$P_{WT}(W)$},ymax =150000,xmax=24,legend pos=north east, xmajorgrids,ymajorgrids]
\addplot [smooth,mark=o,color=green!80!blue, thick]
coordinates{
    (0,130000)(1,130000)(2,130000)(3,130000)(4,130000)(5,130000)(6,94770)(7,94770)(8,28080)(9,66560)(10,16250)(11,94770)(12,130000)(13,130000)(14,130000)(15,130000)(16,130000)(17,130000)(18,130000)(19,130000)(20,130000)(21,130000)(22,130000)(23,130000)
    };    
%\addlegendentry{$P_{WT}max(t)$}
\label{tikz:pwmax}
\addplot[smooth,mark=x, color=blue!80!green, thick] coordinates {(0,129999.626)(1,129999.626)(2,129999.7395)(3,129999.7395)(4,129999.626)(5,129999.4928)(6,94769.9999)(7,94769.9999)(8,28079.9999)(9,66559.9999)(10,16249.9998)(11,94769.9999)(12,129997.7414)(13,129997.7418)(14,129997.7418)(15,129997.7418)(16,129998.0251)(17,129998.0251)(18,129997.1217)(19,129996.0628)(20,129997.1224)(21,129998.2905)(22,129998.7664)(23,129999.4928)};
\label{tikz:pw}
%\addlegendentry{$P_{WT}(t)$}

\end{axis}
\end{tikzpicture}

%% file: Figures/Baterias_PowSale.tex
\begin{tikzpicture}
\begin{axis}[width=9cm,height=5cm,xlabel={Time (h)},ylabel={Energy (Wh)},xmax=24,ymax =45000,legend pos=north west, xmajorgrids,ymajorgrids,legend style={nodes={scale=0.6, transform shape},legend columns=3}]
%\addplot [smooth,mark=*,color=blue!30!green,very thick]
%coordinates{    (0,56070)(1,56070)(2,46725)(3,46725)(4,56070)(5,65415)(6,84105)(7,102795)(8,121485)(9,130830)(10,140175)(11,149520)(12,140175)(13,140175)(14,140175)(15,140175)(16,130830)(17,130830)(18,158865)(19,186900)(20,158865)(21,121485)(22,102795)(23,65415)};    
%\addlegendentry{$Demand$}
\addplot [mark=10-pointed star, color=blue!30!green,thick]coordinates{(0,20442.4758)(1,21885.2711)(2,22581.2659)(3,25534.0725)(4,26977.0129)(5,29165.2646)(6,32841.1375)(7,37999.9987)(8,37900.0349)(9,37800.2693)(10,37999.9991)(11,28882.2264)(12,19000.0012)(13,23092.1579)(14,27184.3666)(15,31276.5861)(16,34638.3103)(17,37999.9989)(18,27363.3077)(19,19000.0012)(20,19200.0283)(21,19100.0485)(22,19000.0016)(23,19000.0001)};
\addlegendentry{$E_\text{Node11}$}

\addplot[mark=triangle*, color=blue!60,thick] coordinates {(0,20469.1635)(1,21938.4562)(2,22686.3654)(3,25691.0332)(4,27160.4533)(5,29348.2503)(6,32965.1575)(7,37999.9988)(8,37900.0906)(9,37800.5617)(10,37999.9991)(11,28860.9565)(12,19000.001)(13,23076.4987)(14,27153.0785)(15,31229.6834)(16,34614.8684)(17,37999.9988)(18,27433.7864)(19,19000.0012)(20,19200.0752)(21,19100.0662)(22,19000.0019)(23,19000.0001)};
\addlegendentry{$E_\text{Node13}$}

\addplot[mark=square*, color=orange,thick] coordinates {(0,19862.9485)(1,20726.1208)(2,20626.5968)(3,22385.4521)(4,23248.9647)(5,25464.5479)(6,30384.7502)(7,37999.9985)(8,37900.0156)(9,37800.199)(10,37999.999)(11,29189.1723)(12,19000.0012)(13,23327.1661)(14,27654.4036)(15,31981.6583)(16,34990.8596)(17,37999.9987)(18,26455.5066)(19,19000.0014)(20,19200.2249)(21,19100.1723)(22,19000.0023)(23,19000.0001)};

\addlegendentry{$E_\text{Node14}$}
\addplot[mark=diamond*, color=cyan,thick] coordinates {(0,19482.0635)(1,19964.9483)(2,19865.1091)(3,20913.5926)(4,21397.7799)(5,23562.2744)(6,29099.1693)(7,37999.9984)(8,37899.9871)(9,37800.1144)(10,37999.999)(11,29359.9507)(12,19000.0013)(13,23457.5499)(14,27915.1784)(15,32372.8238)(16,35186.451)(17,37999.9986)(18,25946.3482)(19,19000.0014)(20,19200.2617)(21,19100.1773)(22,19000.0027)(23,19000.0001)};
\addlegendentry{$E_\text{Node18}$}

\addplot[mark=pentagon*, color=blue!50!green,thick] coordinates {(0,20694.7769)(1,22389.6912)(2,23623.3441)(3,27116.0113)(4,28811.0175)(5,30965.0855)(6,34030.9605)(7,37999.9989)(8,37900.2262)(9,37801.5942)(10,37999.9991)(11,28729.9224)(12,19000.001)(13,22976.1667)(14,26952.4142)(15,30928.6808)(16,34464.3586)(17,37999.9989)(18,27821.4842)(19,19000.0011)(20,19199.9392)(21,19099.9997)(22,19000.0016)(23,19000.0001)};
\addlegendentry{$E_\text{Node19}$}

%\addplot[smooth,mark=*, color=green!50!blue, thick] coordinates {(0,115777.1837)(1,121901.7824)(2,123387.9288)(3,124872.4241)(4,125206.8091)(5,124713.4418)(6,124224.5764)(7,124646.095)(8,126950.6893)(9,130134.1242)(10,134081.5109)(11,137131.5921)(12,140295.254)(13,142941.1003)(14,145574.7786)(15,146966.2779)(16,146471.4127)(17,145969.7232)(18,140293.3686)(19,129084.779)(20,123406.8133)(21,122900.6321)(22,122397.5906)(23,121900.739)};

%\addlegendentry{$EB_{tot}$}

\addplot[smooth,dashed, color=green] coordinates {(0,38000)(1,38000)(2,38000)(3,38000)(4,38000)(5,38000)(6,38000)(7,38000)(8,38000)(9,38000)(10,38000)(11,38000)(12,38000)(13,38000)(14,38000)(15,38000)(16,38000)(17,38000)(18,38000)(19,38000)(20,38000)(21,38000)(22,38000)(23,38000)};

\addlegendentry{$E_\text{max-min}$}
\addplot[smooth,dashed, color=green] coordinates {(0,19000)(1,19000)(2,19000)(3,19000)(4,19000)(5,19000)(6,19000)(7,19000)(8,19000)(9,19000)(10,19000)(11,19000)(12,19000)(13,19000)(14,19000)(15,19000)(16,19000)(17,19000)(18,19000)(19,19000)(20,19000)(21,19000)(22,19000)(23,19000)};

%\addlegendentry{$E_{max}$}

\end{axis}
\end{tikzpicture}

%% file: Figures/GenvsDem.tex
\begin{tikzpicture}
\begin{axis}[width=9cm, height=5cm, xlabel={Time (h)},ylabel={ Power(kW)},xmax=24, xmajorgrids,ymajorgrids]
\addplot [smooth,mark=o,color=green!80!blue, thick]
coordinates{
    (0,56.07)(1,56.07)(2,46.725)(3,46.725)(4,56.07)(5,65.415)(6,84.105)(7,102.795)(8,121.485)(9,130.83)(10,140.175)(11,149.52)(12,140.175)(13,140.175)(14,140.175)(15,140.175)(16,130.83)(17,130.83)(18,158.865)(19,186.9)(20,158.865)(21,121.485)(22,102.795)(23,65.415)
    };
\label{tikz:totaldemand}    
%\addlegendentry{$Demand$}
\addplot[smooth,mark=triangle, color=orange, thick] coordinates {(0,55.4956)(1,55.5161)(2,46.3047)(3,46.3092)(4,55.5148)(5,64.7343)(6,83.4473)(7,101.9971)(8,121.1812)(9,130.1117)(10,139.9525)(11,148.3531)(12,138.7858)(13,138.7749)(14,138.8166)(15,138.8356)(16,129.5174)(17,129.3323)(18,156.1361)(19,183.6895)(20,140.1082)(21,119.7714)(22,101.5746)(23,64.7484)};
\label{tikz:generation}
%\addlegendentry{$P_{ER}$}
%\addplot[smooth,mark=x, color=blue!80, very thick] coordinates {(0,-27266.1638)(1,-7239.2603)(2,-2527.6989)(3,-2526.1751)(4,-1357.1919)(5,-517.6162)(6,-523.6542)(7,-1450.1468)(8,-3900.5678)(9,-10401.3616)(10,-15298.3514)(11,-14676.1773)(12,-12841.0742)(13,-10407.9756)(14,-8375.2381)(15,-2438.4689)(16,-517.4869)(17,-508.1953)(18,4585.1899)(19,9998.1331)(20,4586.7587)(21,-503.032)(22,-506.6334)(23,-514.5192)};

%\addlegendentry{$P_{B}$}
\addplot[smooth,mark=square, color=blue!80, thick] coordinates {(0,0)(1,0)(2,0)(3,0)(4,0)(5,0)(6,0)(7,0)(8,0)(9,0)(10,0)(11,0)(12,0)(13,0)(14,0)(15,0)(16,0)(17,0)(18,0.00015416)(19,0.00012463)(20,15.4466)(21,0)(22,0)(23,0)};
\label{tikz:pslack}
%\addlegendentry{$P_{slack}$}
\end{axis}
\end{tikzpicture}

%% file: Figures/SvSmax.tex
\begin{tikzpicture}
\begin{axis}[width=9cm, height=5cm,xlabel={Time (h)},ylabel={$P_{PV}(W)$},xmax=24, xmajorgrids,ymajorgrids]
\addplot [color=green!80!blue,mark=o,thick]
coordinates{
    (0,0)(1,0)(2,0)(3,0)(4,0)(5,0)(6,30000)(7,60000)(8,90000)(9,120000)(10,142500)(11,150000)(12,135000)(13,127500)(14,120000)(15,90000)(16,60000)(17,30000)(18,0)(19,0)(20,0)(21,0)(22,0)(23,0)
    };    
%\addlegendentry{$P_{PV}max(t)$}
\label{tikz:pvmax2}
\addplot[mark=x, color=blue!80!green, thick] coordinates {(0,0)(1,0)(2,0)(3,0)(4,0)(5,0)(6,18646.3791)(7,43350.1075)(8,75521.1877)(9,99508.7888)(10,116509.0385)(11,98778.7954)(12,75036.943)(13,71387.8221)(14,69656.5036)(15,59195.8869)(16,42959.4411)(17,24665.9526)(18,0)(19,0)(20,0)(21,0)(22,0)(23,0)};
\label{tikz:pv2}
%\addlegendentry{$P_{PV}(t)$}
\end{axis}
\end{tikzpicture}

%% file: Figures/WvWmax.tex
\begin{tikzpicture}
\begin{axis}[width=9cm, height=5.0cm,xlabel={Time (h)},ylabel={$P_{WT}(W)$},ymax =140000,xmax=25,legend pos=north east, xmajorgrids,ymajorgrids,legend style={nodes={scale=1, transform shape}}]
\addplot [smooth,mark=o,color=green!80!blue,thick]
coordinates{
    (0,130000)(1,130000)(2,130000)(3,130000)(4,130000)(5,130000)(6,94770)(7,94770)(8,28080)(9,66560)(10,16250)(11,94770)(12,130000)(13,130000)(14,130000)(15,130000)(16,130000)(17,130000)(18,130000)(19,130000)(20,130000)(21,130000)(22,130000)(23,130000)
    };    
%\addlegendentry{$P_{WT}max(t)$}
\label{tikz:pwmax2}
\addplot[smooth,mark=x, color=blue!80!green, thick] coordinates {(0,104800.2823)(1,100704.0138)(2,91275.7779)(3,89907.1222)(4,100506.995)(5,108501.0805)(6,89420.3776)(7,89757.9334)(8,24579.9356)(9,62792.5511)(10,13009.2933)(11,87618.6781)(12,113108.0818)(13,114479.114)(14,114113.5703)(15,117091.8278)(16,121034.9089)(17,125963.6912)(18,129997.1313)(19,129996.0678)(20,129997.1744)(21,128708.2369)(22,124302.9595)(23,106871.8017)};

%\addlegendentry{$P_{WT}(t)$}
\label{tikz:pw2}
\end{axis}
\end{tikzpicture}

%% file: Figures/Baterias.tex
\begin{tikzpicture}
\begin{axis}[width=9cm, height=5 cm, xlabel={Time (h)},ylabel={Energy (Wh)},ymax =45000,xmax=24,legend pos=north west, xmajorgrids,ymajorgrids,legend style={nodes={scale=0.6, transform shape},legend columns=3}]
%\addplot [smooth,mark=*,color=blue!30!green,very thick]
%coordinates{    (0,56070)(1,56070)(2,46725)(3,46725)(4,56070)(5,65415)(6,84105)(7,102795)(8,121485)(9,130830)(10,140175)(11,149520)(12,140175)(13,140175)(14,140175)(15,140175)(16,130830)(17,130830)(18,158865)(19,186900)(20,158865)(21,121485)(22,102795)(23,65415)};    
%\addlegendentry{$Demand$}
\addplot[mark=10-pointed star, color=blue!30!green,thick] coordinates {(0,23898.5442)(1,26539.9906)(2,28579.5941)(3,29934.4745)(4,31554.7053)(5,33832.0651)(6,33389.365)(7,34902.906)(8,26998.519)(9,30341.7415)(10,24358.7635)(11,25772.0892)(12,28611.3212)(13,30656.5153)(14,32648.171)(15,34035.8763)(16,36123.5851)(17,37999.9991)(18,32705.4092)(19,21019.8669)(20,19000.0015)(21,19536.4351)(22,21009.4365)(23,22675.9141)};
\addlegendentry{$EB_{Node11}$}

\addplot[mark=triangle*, color=blue!60,thick] coordinates {(0,23535.6888)(1,26197.4001)(2,28379.5175)(3,30096.019)(4,32051.8661)(5,34273.9013)(6,33686.0277)(7,35133.33)(8,27511.8204)(9,29855.3657)(10,23560.9804)(11,25335.4256)(12,28431.999)(13,30840.2353)(14,32976.2776)(15,34371.004)(16,36301.9333)(17,37999.999)(18,32755.3401)(19,20968.4788)(20,19000.002)(21,19578.333)(22,21124.4634)(23,22766.1467)};

\addlegendentry{$EB_{Node13}$}
\addplot[mark=square*, color=orange,thick] coordinates {(0,23286.7372)(1,25835.9881)(2,28000.1378)(3,29803.3824)(4,31546.91)(5,33145.5101)(6,33275.6865)(7,34976.5783)(8,27424.811)(9,29862.0084)(10,23717.882)(11,25520.0369)(12,28512.4996)(13,30801.293)(14,32912.968)(15,34365.139)(16,36373.0078)(17,37999.999)(18,32105.2227)(19,21622.8172)(20,19000.0031)(21,19791.7746)(22,21135.2519)(23,22276.2629)};

\addlegendentry{$EB_{Node14}$}
\addplot[mark=diamond*, color=cyan, thick] coordinates {(0,23269.6145)(1,25821.5257)(2,27980.846)(3,29802.509)(4,31547.4471)(5,32972.3142)(6,33203.0798)(7,34929.2737)(8,27361.8713)(9,29749.7041)(10,23685.7064)(11,25552.5473)(12,28512.9318)(13,30765.1532)(14,32863.3489)(15,34331.6201)(16,36341.4295)(17,37999.999)(18,31767.0761)(19,21962.0293)(20,19000.0041)(21,19701.9009)(22,21275.7554)(23,22168.3964)};

\addlegendentry{$EB_{Node18}$}
\addplot[mark=pentagon*, color=blue!50!green,thick] coordinates {(0,23614.9022)(1,26303.3268)(2,28472.2755)(3,30047.3616)(4,31982.8432)(5,34245.8629)(6,33645.2582)(7,35095.7458)(8,27417.8202)(9,29853.8183)(10,23601.5868)(11,25315.8428)(12,28424.8272)(13,30833.712)(14,32981.5142)(15,34335.0455)(16,36193.6986)(17,37999.9991)(18,33013.8498)(19,20709.4405)(20,19000.002)(21,19563.3879)(22,20973.9993)(23,22423.974)};

\addlegendentry{$EB_{Node19}$}
%\addplot[smooth,mark=*, color=green!50!blue, thick] coordinates {(0,115777.1837)(1,121901.7824)(2,123387.9288)(3,124872.4241)(4,125206.8091)(5,124713.4418)(6,124224.5764)(7,124646.095)(8,126950.6893)(9,130134.1242)(10,134081.5109)(11,137131.5921)(12,140295.254)(13,142941.1003)(14,145574.7786)(15,146966.2779)(16,146471.4127)(17,145969.7232)(18,140293.3686)(19,129084.779)(20,123406.8133)(21,122900.6321)(22,122397.5906)(23,121900.739)};

%\addlegendentry{$EB_{tot}$}

\addplot[smooth,dashed, color=green] coordinates {(0,38000)(1,38000)(2,38000)(3,38000)(4,38000)(5,38000)(6,38000)(7,38000)(8,38000)(9,38000)(10,38000)(11,38000)(12,38000)(13,38000)(14,38000)(15,38000)(16,38000)(17,38000)(18,38000)(19,38000)(20,38000)(21,38000)(22,38000)(23,38000)};

\addlegendentry{$E_\text{max-min}$}
\addplot[smooth,dashed, color=green] coordinates {(0,19000)(1,19000)(2,19000)(3,19000)(4,19000)(5,19000)(6,19000)(7,19000)(8,19000)(9,19000)(10,19000)(11,19000)(12,19000)(13,19000)(14,19000)(15,19000)(16,19000)(17,19000)(18,19000)(19,19000)(20,19000)(21,19000)(22,19000)(23,19000)};

%\addlegendentry{$E_{max}$}

\end{axis}
\end{tikzpicture}

%% file: Figures/IEEE123.tex
\scriptsize
\centering
\begin{tikzpicture}[x=0.20mm, y = 0.20mm]
 \fill (99,-323) circle[radius=2pt] node [inner sep = 0pt](N1){}; 
 \fill (93,-284) circle[radius=2pt] node [inner sep = 0pt](N2){}; 
 \fill (99,-383) circle[radius=2pt] node [inner sep = 0pt](N3){}; 
 \fill (99,-402) circle[radius=2pt] node [inner sep = 0pt](N4){}; 
 \fill (126,-383) circle[radius=2pt] node [inner sep = 0pt](N5){}; 
 \fill (157,-384) circle[radius=2pt] node [inner sep = 0pt](N6){}; 
 \fill (123,-318) circle[radius=2pt] node [inner sep = 0pt](N7){}; 
 \fill (155,-313) circle[radius=2pt] node [inner sep = 0pt](N8){}; 
 \fill (146,-275) circle[radius=2pt] node [inner sep = 0pt](N9){}; 
 \fill (110,-272) circle[radius=2pt] node [inner sep = 0pt](N10){}; 
 \fill (64,-254) circle[radius=2pt] node [inner sep = 0pt](N11){}; 
 \fill (140,-337) circle[radius=2pt] node [inner sep = 0pt](N12){}; 
 \fill (187,-307) circle[radius=2pt] node [inner sep = 0pt](N13){}; 
 \fill (124,-249) circle[radius=2pt] node [inner sep = 0pt](N14){}; 
 \fill (207,-369) circle[radius=2pt] node [inner sep = 0pt](N15){}; 
 \fill (219,-397) circle[radius=2pt] node [inner sep = 0pt](N16){}; 
 \fill (234,-360) circle[radius=2pt] node [inner sep = 0pt](N17){}; 
 \fill (149,-184) circle[radius=2pt] node [inner sep = 0pt](N18){}; 
 \fill (111,-193) circle[radius=2pt] node [inner sep = 0pt](N19){}; 
 \fill (71,-206) circle[radius=2pt] node [inner sep = 0pt](N20){}; 
 \fill (137,-146) circle[radius=2pt] node [inner sep = 0pt](N21){}; 
 \fill (63,-165) circle[radius=2pt] node [inner sep = 0pt](N22){}; 
 \fill (123,-101) circle[radius=2pt] node [inner sep = 0pt](N23){}; 
 \fill (65,-116) circle[radius=2pt] node [inner sep = 0pt](N24){}; 
 \fill (114,-70) circle[radius=2pt] node [inner sep = 0pt](N25){}; 
 \fill (76,-78) circle[radius=2pt] node [inner sep = 0pt](N26){}; 
 \fill (31,-87) circle[radius=2pt] node [inner sep = 0pt](N27){}; 
 \fill (107,-44) circle[radius=2pt] node [inner sep = 0pt](N28){}; 
 \fill (99,-17) circle[radius=2pt] node [inner sep = 0pt](N29){}; 
 \fill (130,-14) circle[radius=2pt] node [inner sep = 0pt](N30){}; 
 \fill (64,-44) circle[radius=2pt] node [inner sep = 0pt](N31){}; 
 \fill (55,-17) circle[radius=2pt] node [inner sep = 0pt](N32){}; 
 \fill (22,-42) circle[radius=2pt] node [inner sep = 0pt](N33){}; 
 \fill (197,-337) circle[radius=2pt] node [inner sep = 0pt](N34){}; 
 \fill (204,-171) circle[radius=2pt] node [inner sep = 0pt](N35){}; 
 \fill (260,-188) circle[radius=2pt] node [inner sep = 0pt](N36){}; 
 \fill (207,-202) circle[radius=2pt] node [inner sep = 0pt](N37){}; 
 \fill (294,-176) circle[radius=2pt] node [inner sep = 0pt](N38){}; 
 \fill (322,-168) circle[radius=2pt] node [inner sep = 0pt](N39){}; 
 \fill (195,-148) circle[radius=2pt] node [inner sep = 0pt](N40){}; 
 \fill (244,-133) circle[radius=2pt] node [inner sep = 0pt](N41){}; 
 \fill (190,-124) circle[radius=2pt] node [inner sep = 0pt](N42){}; 
 \fill (245,-104) circle[radius=2pt] node [inner sep = 0pt](N43){}; 
 \fill (181,-98) circle[radius=2pt] node [inner sep = 0pt](N44){}; 
 \fill (218,-88) circle[radius=2pt] node [inner sep = 0pt](N45){}; 
 \fill (258,-75) circle[radius=2pt] node [inner sep = 0pt](N46){}; 
 \fill (170,-66) circle[radius=2pt] node [inner sep = 0pt](N47){}; 
 \fill (141,-74) circle[radius=2pt] node [inner sep = 0pt](N48){}; 
 \fill (209,-54) circle[radius=2pt] node [inner sep = 0pt](N49){}; 
 \fill (245,-44) circle[radius=2pt] node [inner sep = 0pt](N50){}; 
 \fill (276,-36) circle[radius=2pt] node [inner sep = 0pt](N51){}; 
 \fill (268,-293) circle[radius=2pt] node [inner sep = 0pt](N52){}; 
 \fill (298,-288) circle[radius=2pt] node [inner sep = 0pt](N53){}; 
 \fill (323,-282) circle[radius=2pt] node [inner sep = 0pt](N54){}; 
 \fill (350,-281) circle[radius=2pt] node [inner sep = 0pt](N55){}; 
 \fill (377,-276) circle[radius=2pt] node [inner sep = 0pt](N56){}; 
 \fill (313,-238) circle[radius=2pt] node [inner sep = 0pt](N57){}; 
 \fill (281,-247) circle[radius=2pt] node [inner sep = 0pt](N58){}; 
 \fill (245,-255) circle[radius=2pt] node [inner sep = 0pt](N59){}; 
 \fill (381,-224) circle[radius=2pt] node [inner sep = 0pt](N60){}; 
 \fill (394,-256) circle[radius=2pt] node [inner sep = 0pt](N61){}; 
 \fill (380,-154) circle[radius=2pt] node [inner sep = 0pt](N62){}; 
 \fill (369,-119) circle[radius=2pt] node [inner sep = 0pt](N63){}; 
 \fill (358,-90) circle[radius=2pt] node [inner sep = 0pt](N64){}; 
 \fill (318,-107) circle[radius=2pt] node [inner sep = 0pt](N65){}; 
 \fill (326,-133) circle[radius=2pt] node [inner sep = 0pt](N66){}; 
 \fill (462,-213) circle[radius=2pt] node [inner sep = 0pt](N67){}; 
 \fill (485,-200) circle[radius=2pt] node [inner sep = 0pt](N68){}; 
 \fill (510,-186) circle[radius=2pt] node [inner sep = 0pt](N69){}; 
 \fill (534,-171) circle[radius=2pt] node [inner sep = 0pt](N70){}; 
 \fill (560,-152) circle[radius=2pt] node [inner sep = 0pt](N71){}; 
 \fill (475,-250) circle[radius=2pt] node [inner sep = 0pt](N72){}; 
 \fill (503,-236) circle[radius=2pt] node [inner sep = 0pt](N73){}; 
 \fill (527,-220) circle[radius=2pt] node [inner sep = 0pt](N74){}; 
 \fill (552,-206) circle[radius=2pt] node [inner sep = 0pt](N75){}; 
 \fill (487,-293) circle[radius=2pt] node [inner sep = 0pt](N76){}; 
 \fill (505,-285) circle[radius=2pt] node [inner sep = 0pt](N77){}; 
 \fill (522,-277) circle[radius=2pt] node [inner sep = 0pt](N78){}; 
 \fill (540,-266) circle[radius=2pt] node [inner sep = 0pt](N79){}; 
 \fill (524,-308) circle[radius=2pt] node [inner sep = 0pt](N80){}; 
 \fill (527,-343) circle[radius=2pt] node [inner sep = 0pt](N81){}; 
 \fill (528,-383) circle[radius=2pt] node [inner sep = 0pt](N82){}; 
 \fill (560,-383) circle[radius=2pt] node [inner sep = 0pt](N83){}; 
 \fill (563,-318) circle[radius=2pt] node [inner sep = 0pt](N84){}; 
 \fill (561,-248) circle[radius=2pt] node [inner sep = 0pt](N85){}; 
 \fill (492,-361) circle[radius=2pt] node [inner sep = 0pt](N86){}; 
 \fill (437,-365) circle[radius=2pt] node [inner sep = 0pt](N87){}; 
 \fill (429,-324) circle[radius=2pt] node [inner sep = 0pt](N88){}; 
 \fill (395,-369) circle[radius=2pt] node [inner sep = 0pt](N89){}; 
 \fill (391,-332) circle[radius=2pt] node [inner sep = 0pt](N90){}; 
 \fill (352,-375) circle[radius=2pt] node [inner sep = 0pt](N91){}; 
 \fill (349,-341) circle[radius=2pt] node [inner sep = 0pt](N92){}; 
 \fill (312,-376) circle[radius=2pt] node [inner sep = 0pt](N93){}; 
 \fill (304,-326) circle[radius=2pt] node [inner sep = 0pt](N94){}; 
 \fill (280,-381) circle[radius=2pt] node [inner sep = 0pt](N95){}; 
 \fill (274,-333) circle[radius=2pt] node [inner sep = 0pt](N96){}; 
 \fill (447,-174) circle[radius=2pt] node [inner sep = 0pt](N97){}; 
 \fill (467,-160) circle[radius=2pt] node [inner sep = 0pt](N98){}; 
 \fill (487,-148) circle[radius=2pt] node [inner sep = 0pt](N99){}; 
 \fill (522,-122) circle[radius=2pt] node [inner sep = 0pt](N100){}; 
 \fill (434,-136) circle[radius=2pt] node [inner sep = 0pt](N101){}; 
 \fill (459,-119) circle[radius=2pt] node [inner sep = 0pt](N102){}; 
 \fill (489,-101) circle[radius=2pt] node [inner sep = 0pt](N103){}; 
 \fill (516,-82) circle[radius=2pt] node [inner sep = 0pt](N104){}; 
 \fill (424,-107) circle[radius=2pt] node [inner sep = 0pt](N105){}; 
 \fill (451,-89) circle[radius=2pt] node [inner sep = 0pt](N106){}; 
 \fill (492,-66) circle[radius=2pt] node [inner sep = 0pt](N107){}; 
 \fill (413,-76) circle[radius=2pt] node [inner sep = 0pt](N108){}; 
 \fill (445,-61) circle[radius=2pt] node [inner sep = 0pt](N109){}; 
 \fill (476,-37) circle[radius=2pt] node [inner sep = 0pt](N110){}; 
 \fill (429,-38) circle[radius=2pt] node [inner sep = 0pt](N111){}; 
 \fill (505,-38) circle[radius=2pt] node [inner sep = 0pt](N112){}; 
 \fill (540,-37) circle[radius=2pt] node [inner sep = 0pt](N113){}; 
 \fill (572,-37) circle[radius=2pt] node [inner sep = 0pt](N114){}; 
 \fill (70,-323) circle[radius=2pt] node [inner sep = 0pt](N116){}; 
 \fill (30,-323) circle[radius=2pt] node [inner sep = 0pt](N117){}; 
 \fill (420,-220) circle[radius=2pt] node [inner sep = 0pt](N120){}; 
 \fill (180,-175) circle[radius=2pt] node [inner sep = 0pt](N115){}; 
 \fill (220,-300) circle[radius=2pt] node [inner sep = 0pt](N119){}; 
 \fill (420,-256) circle[radius=2pt] node [inner sep = 0pt](N129){}; 
 \fill (441,-155) circle[radius=2pt] node [inner sep = 0pt](N122){}; 
 \fill (310,-36) circle[radius=2pt] node [inner sep = 0pt](N118){}; 
 \fill (150,-14) circle[radius=2pt] node [inner sep = 0pt](N123){}; 
 \fill (360,-36) circle[radius=2pt] node [inner sep = 0pt](N125){}; 
 \fill (550,-122) circle[radius=2pt] node [inner sep = 0pt](N127){}; 
 \node at (N1)[yshift=-6,xshift=-6] {1}; 
 \node at (N2)[yshift=6] {2}; 
 \node at (N3)[xshift=-6] {3}; 
 \node at (N4)[yshift=-6] {4}; 
 \node at (N5)[yshift=6] {5}; 
 \node at (N6)[yshift=6] {6}; 
 \node at (N7)[yshift=6] {7}; 
 \node at (N8)[yshift=-6] {8}; 
 \node at (N9)[yshift=6] {9}; 
 \node at (N10)[yshift=-6] {10}; 
 \node at (N11)[yshift=6] {11}; 
 \node at (N12)[yshift=-6] {12}; 
 \node at (N13)[yshift=6,xshift=6] {13}; 
 \node at (N14)[yshift=6] {14}; 
 \node at (N15)[xshift=-6] {15}; 
 \node at (N16)[yshift=-6] {16}; 
 \node at (N17)[yshift=6] {17}; 
 \node at (N18)[yshift=6,xshift=4] {18}; 
 \node at (N19)[yshift=6] {19}; 
 \node at (N20)[yshift=-6] {20}; 
 \node at (N21)[xshift=6] {21}; 
 \node at (N22)[yshift=6] {22}; 
 \node at (N23)[xshift=6] {23}; 
 \node at (N24)[yshift=6] {24}; 
 \node at (N25)[xshift=6] {25}; 
 \node at (N26)[yshift=-6] {26}; 
 \node at (N27)[yshift=-6] {27}; 
 \node at (N28)[xshift=6] {28}; 
 \node at (N29)[yshift=6] {29}; 
 \node at (N30)[yshift=-6] {30}; 
 \node at (N31)[xshift=6] {31}; 
 \node at (N32)[yshift=6] {32}; 
 \node at (N33)[yshift=6] {33}; 
 \node at (N34)[xshift=6] {34}; 
 \node at (N35)[xshift=6,yshift=5] {35}; 
 \node at (N36)[yshift=6] {36}; 
 \node at (N37)[yshift=6] {37}; 
 \node at (N38)[yshift=6] {38}; 
 \node at (N39)[yshift=-6] {39}; 
 \node at (N40)[xshift=-6] {40}; 
 \node at (N41)[yshift=6] {41}; 
 \node at (N42)[xshift=-6] {42}; 
 \node at (N43)[yshift=6] {43}; 
 \node at (N44)[xshift=-6] {44}; 
 \node at (N45)[yshift=6] {45}; 
 \node at (N46)[yshift=6] {46}; 
 \node at (N47)[yshift=6] {47}; 
 \node at (N48)[yshift=6] {48}; 
 \node at (N49)[yshift=6] {49}; 
 \node at (N50)[yshift=-6] {50}; 
 \node at (N51)[yshift=6] {51}; 
 \node at (N52)[yshift=6] {52}; 
 \node at (N53)[yshift=6] {53}; 
 \node at (N54)[yshift=-6] {54}; 
 \node at (N55)[yshift=6] {55}; 
 \node at (N56)[yshift=6] {56}; 
 \node at (N57)[yshift=6] {57}; 
 \node at (N58)[yshift=6] {58}; 
 \node at (N59)[yshift=6] {59}; 
 \node at (N60)[yshift=6,xshift=6] {60}; 
 \node at (N61)[yshift=-6] {61}; 
 \node at (N62)[xshift=6] {62}; 
 \node at (N63)[xshift=6] {63}; 
 \node at (N64)[yshift=6] {64}; 
 \node at (N65)[yshift=6] {65}; 
 \node at (N66)[yshift=-6] {66}; 
 \node at (N67)[yshift=6,xshift=-8] {67}; 
 \node at (N68)[yshift=6] {68}; 
 \node at (N69)[yshift=6] {69}; 
 \node at (N70)[yshift=6] {70}; 
 \node at (N71)[yshift=6] {71}; 
 \node at (N72)[xshift=-6] {72}; 
 \node at (N73)[yshift=6] {73}; 
 \node at (N74)[yshift=6] {74}; 
 \node at (N75)[yshift=6] {75}; 
 \node at (N76)[xshift=-6] {76}; 
 \node at (N77)[yshift=6] {77}; 
 \node at (N78)[yshift=6] {78}; 
 \node at (N79)[yshift=6] {79}; 
 \node at (N80)[xshift=6] {80}; 
 \node at (N81)[xshift=-6] {81}; 
 \node at (N82)[xshift=-6] {82}; 
 \node at (N83)[yshift=6] {83}; 
 \node at (N84)[yshift=-6] {84}; 
 \node at (N85)[yshift=6] {85}; 
 \node at (N86)[yshift=6, xshift=-7] {86}; 
 \node at (N87)[yshift=-6] {87}; 
 \node at (N88)[yshift=6] {88}; 
 \node at (N89)[yshift=-6] {89}; 
 \node at (N90)[yshift=6] {90}; 
 \node at (N91)[yshift=-6] {91}; 
 \node at (N92)[yshift=6] {92}; 
 \node at (N93)[yshift=-6] {93}; 
 \node at (N94)[yshift=6] {94}; 
 \node at (N95)[yshift=-6] {95}; 
 \node at (N96)[yshift=6] {96}; 
 \node at (N97)[xshift=-7] {97}; 
 \node at (N98)[yshift=6] {98}; 
 \node at (N99)[yshift=6] {99}; 
 \node at (N100)[yshift=6] {100}; 
 \node at (N101)[xshift=-9] {101}; 
 \node at (N102)[yshift=6] {102}; 
 \node at (N103)[yshift=6] {103}; 
 \node at (N104)[yshift=6] {104}; 
 \node at (N105)[xshift=-9] {105}; 
 \node at (N106)[yshift=6] {106}; 
 \node at (N107)[yshift=6] {107}; 
 \node at (N108)[yshift=6] {108}; 
 \node at (N109)[yshift=6] {109}; 
 \node at (N110)[yshift=6] {110}; 
 \node at (N111)[yshift=6] {111}; 
 \node at (N112)[yshift=6] {112}; 
 \node at (N113)[yshift=6] {113}; 
 \node at (N114)[yshift=6] {114}; 
 \node at (N116)[yshift=6] {149}; 
 \node at (N117)[yshift=6] {150}; 
 \node at (N120)[yshift=6] {160}; 
 \node at (N115)[yshift=6] {135}; 
 \node at (N119)[yshift=6] {152}; 
 \node at (N129)[yshift=6] {610}; 
 \node at (N122)[xshift=-9] {197}; 
 \node at (N118)[yshift=6] {151}; 
 \node at (N123)[yshift=6] {250}; 
 \node at (N125)[yshift=6] {300}; 
 \node at (N127)[yshift=6] {450};

 \draw[thick,black] (N116) -- (N1); 
 \draw[thick,gray] (N1) -- (N2); 
 \draw[thick,gray] (N1) -- (N3); 
 \draw[thick,gray] (N3) -- (N4); 
 \draw[thick,gray] (N3) -- (N5); 
 \draw[thick,gray] (N5) -- (N6); 
 \draw[thick,black] (N1) -- (N7); 
 \draw[thick,black] (N7) -- (N8); 
 \draw[thick,gray] (N8) -- (N9); 
 \draw[thick,gray] (N14) -- (N10); 
 \draw[thick,gray] (N14) -- (N11); 
 \draw[thick,gray] (N8) -- (N12); 
 \draw[thick,black] (N8) -- (N13); 
 \draw[thick,gray] (N9) -- (N14); 
 \draw[thick,gray] (N34) -- (N15); 
 \draw[thick,gray] (N15) -- (N16); 
 \draw[thick,gray] (N15) -- (N17); 
 \draw[thick,black] (N13) -- (N18); 
 \draw[thick,gray] (N18) -- (N19); 
 \draw[thick,gray] (N19) -- (N20); 
 \draw[thick,black] (N18) -- (N21); 
 \draw[thick,gray] (N21) -- (N22); 
 \draw[thick,black] (N21) -- (N23); 
 \draw[thick,gray] (N23) -- (N24); 
 \draw[thick,black] (N23) -- (N25); 
 \draw[thick,gray] (N25) -- (N26); 
 \draw[thick,gray] (N26) -- (N27); 
 \draw[thick,black] (N25) -- (N28); 
 \draw[thick,black] (N28) -- (N29); 
 \draw[thick,black] (N29) -- (N30); 
 \draw[thick,gray] (N26) -- (N31); 
 \draw[thick,gray] (N31) -- (N32); 
 \draw[thick,gray] (N27) -- (N33); 
 \draw[thick,gray] (N13) -- (N34); 
 \draw[thick,black] (N115) -- (N35); 
 \draw[thick,gray] (N35) -- (N36); 
 \draw[thick,gray] (N36) -- (N37); 
 \draw[thick,gray] (N36) -- (N38); 
 \draw[thick,gray] (N38) -- (N39); 
 \draw[thick,black] (N35) -- (N40); 
 \draw[thick,gray] (N40) -- (N41); 
 \draw[thick,black] (N40) -- (N42); 
 \draw[thick,gray] (N42) -- (N43); 
 \draw[thick,black] (N42) -- (N44); 
 \draw[thick,gray] (N44) -- (N45); 
 \draw[thick,gray] (N45) -- (N46); 
 \draw[thick,black] (N44) -- (N47); 
 \draw[thick,black] (N47) -- (N48); 
 \draw[thick,black] (N47) -- (N49); 
 \draw[thick,black] (N49) -- (N50); 
 \draw[thick,black] (N50) -- (N51); 
 \draw[thick,black] (N119) -- (N52); 
 \draw[thick,black] (N52) -- (N53); 
 \draw[thick,black] (N53) -- (N54); 
 \draw[thick,black] (N54) -- (N55); 
 \draw[thick,black] (N55) -- (N56); 
 \draw[thick,black] (N54) -- (N57); 
 \draw[thick,gray] (N57) -- (N58); 
 \draw[thick,gray] (N58) -- (N59); 
 \draw[thick,black] (N57) -- (N60); 
 \draw[thick,black] (N60) -- (N61); 
 \draw[thick,black] (N60) -- (N62); 
 \draw[thick,black] (N62) -- (N63); 
 \draw[thick,black] (N63) -- (N64); 
 \draw[thick,black] (N64) -- (N65); 
 \draw[thick,black] (N65) -- (N66); 
 \draw[thick,black] (N120) -- (N67); 
 \draw[thick,gray] (N67) -- (N68); 
 \draw[thick,gray] (N68) -- (N69); 
 \draw[thick,gray] (N69) -- (N70); 
 \draw[thick,gray] (N70) -- (N71); 
 \draw[thick,black] (N67) -- (N72); 
 \draw[thick,gray] (N72) -- (N73); 
 \draw[thick,gray] (N73) -- (N74); 
 \draw[thick,gray] (N74) -- (N75); 
 \draw[thick,black] (N72) -- (N76); 
 \draw[thick,black] (N76) -- (N77); 
 \draw[thick,black] (N77) -- (N78); 
 \draw[thick,black] (N78) -- (N79); 
 \draw[thick,black] (N78) -- (N80); 
 \draw[thick,black] (N80) -- (N81); 
 \draw[thick,black] (N81) -- (N82); 
 \draw[thick,black] (N82) -- (N83); 
 \draw[thick,gray] (N81) -- (N84); 
 \draw[thick,gray] (N84) -- (N85); 
 \draw[thick,black] (N76) -- (N86); 
 \draw[thick,black] (N86) -- (N87); 
 \draw[thick,gray] (N87) -- (N88); 
 \draw[thick,black] (N87) -- (N89); 
 \draw[thick,gray] (N89) -- (N90); 
 \draw[thick,black] (N89) -- (N91); 
 \draw[thick,gray] (N91) -- (N92); 
 \draw[thick,black] (N91) -- (N93); 
 \draw[thick,gray] (N93) -- (N94); 
 \draw[thick,black] (N93) -- (N95); 
 \draw[thick,gray] (N95) -- (N96); 
 \draw[thick,black] (N67) -- (N97); 
 \draw[thick,black] (N97) -- (N98); 
 \draw[thick,black] (N98) -- (N99); 
 \draw[thick,black] (N99) -- (N100); 
 \draw[thick,black] (N122) -- (N101); 
 \draw[thick,gray] (N101) -- (N102); 
 \draw[thick,gray] (N102) -- (N103); 
 \draw[thick,gray] (N103) -- (N104); 
 \draw[thick,black] (N101) -- (N105); 
 \draw[thick,gray] (N105) -- (N106); 
 \draw[thick,gray] (N106) -- (N107); 
 \draw[thick,black] (N105) -- (N108); 
 \draw[thick,gray] (N108) -- (N109); 
 \draw[thick,gray] (N109) -- (N110); 
 \draw[thick,gray] (N110) -- (N111); 
 \draw[thick,gray] (N110) -- (N112); 
 \draw[thick,gray] (N112) -- (N113); 
 \draw[thick,gray] (N113) -- (N114); 
 \draw[thick,black] (N51) -- (N118); 
 \draw[thick,black] (N30) -- (N123); 
 \draw[thick,black] (N108) -- (N125); 
 \draw[thick,black] (N100) -- (N127); 
 \draw[thick,gray] (N61) -- (N129); 
 \draw[thick,black] (N18) -- (N115); 
 \draw[thick,black] (N117) -- (N116); 
 \draw[thick,black] (N13) -- (N119); 
 \draw[thick,black] (N60) -- (N120); 
 \draw[thick,black] (N97) -- (N122); 

% Solar N19
\coordinate[shift={(12,-10)}] (solar) at (N19);
\draw[blue!50!green] (solar) -- (N19);
\draw[blue!50!green, fill=white] (solar) rectangle +(18,-25);
\draw[blue!50!green] (solar) -- +(9,-10) -- +(18,0);

% Solar N20
\coordinate[shift={(-30,20)}] (solar) at (N20);
\draw[blue!50!green] (solar) -- (N20);
\draw[blue!50!green, fill=white] (solar) rectangle +(18,-25);
\draw[blue!50!green] (solar) -- +(9,-10) -- +(18,0);

% Solar N37
\coordinate[shift={(-10,-20)}] (solar) at (N37);
\draw[blue!50!green] (solar) -- (N37);
\draw[blue!50!green, fill=white] (solar) rectangle +(18,-25);
\draw[blue!50!green] (solar) -- +(9,-10) -- +(18,0);

% Solar N43
\coordinate[shift={(30,10)}] (solar) at (N43);
\draw[blue!50!green] (solar) -- (N43);
\draw[blue!50!green, fill=white] (solar) rectangle +(18,-25);
\draw[blue!50!green] (solar) -- +(9,-10) -- +(18,0);

% Solar N50
\coordinate[shift={(-30,30)}] (solar) at (N50);
\draw[blue!50!green] (solar) -- (N50);
\draw[blue!50!green, fill=white] (solar) rectangle +(18,-25);
\draw[blue!50!green] (solar) -- +(9,-10) -- +(18,0);

% Solar N74
\coordinate[shift={(-30,25)}] (solar) at (N74);
\draw[blue!50!green] (solar) -- (N74);
\draw[blue!50!green, fill=white] (solar) rectangle +(18,-25);
\draw[blue!50!green] (solar) -- +(9,-10) -- +(18,0);

% Solar N75
\coordinate[shift={(30,20)}] (solar) at (N75);
\draw[blue!50!green] (solar) -- (N75);
\draw[blue!50!green, fill=white] (solar) rectangle +(18,-25);
\draw[blue!50!green] (solar) -- +(9,-10) -- +(18,0);

% Leyenda Solar
\coordinate (solar) at (50,-430);
\draw[blue!50!green, fill=white] (solar) rectangle +(18,-25);
\draw[blue!50!green] (solar) -- +(9,-10) -- +(18,0);
\node[blue!50!green] at (90,-450) {Solar};

% Bateria N19
\coordinate[shift={(-30,25)}] (bateria) at (N19);
\draw[blue!50!green] (bateria) -- (N19);
\draw[blue!50!green, fill=white] (bateria) rectangle +(14,-25);
\fill[blue!50!green] (bateria)+(1,-8) rectangle +(13,-15);
\fill[blue!50!green] (bateria)+(1,-24) rectangle +(13,-16);

% Bateria N33
\coordinate[shift={(-30,25)}] (bateria) at (N33);
\draw[blue!50!green] (bateria) -- (N33);
\draw[blue!50!green, fill=white] (bateria) rectangle +(14,-25);
\fill[blue!50!green] (bateria)+(1,-8) rectangle +(13,-15);
\fill[blue!50!green] (bateria)+(1,-24) rectangle +(13,-16);

% Bateria N34
\coordinate[shift={(-30,10)}] (bateria) at (N34);
\draw[blue!50!green] (bateria) -- (N34);
\draw[blue!50!green, fill=white] (bateria) rectangle +(14,-25);
\fill[blue!50!green] (bateria)+(1,-8) rectangle +(13,-15);
\fill[blue!50!green] (bateria)+(1,-24) rectangle +(13,-16);

% Bateria N64
\coordinate[shift={(-30,25)}] (bateria) at (N64);
\draw[blue!50!green] (bateria) -- (N64);
\draw[blue!50!green, fill=white] (bateria) rectangle +(14,-25);
\fill[blue!50!green] (bateria)+(1,-8) rectangle +(13,-15);
\fill[blue!50!green] (bateria)+(1,-24) rectangle +(13,-16);

% Bateria N75
\coordinate[shift={(-30,27)}] (bateria) at (N75);
\draw[blue!50!green] (bateria) -- (N75);
\draw[blue!50!green, fill=white] (bateria) rectangle +(14,-25);
\fill[blue!50!green] (bateria)+(1,-8) rectangle +(13,-15);
\fill[blue!50!green] (bateria)+(1,-24) rectangle +(13,-16);

% Leyenda Bateria 
\coordinate (bateria) at (150,-430);
\draw[blue!50!green, fill=white] (bateria) rectangle +(14,-25);
\fill[blue!50!green] (bateria)+(1,-8) rectangle +(13,-15);
\fill[blue!50!green] (bateria)+(1,-24) rectangle +(13,-16);
\node[blue!50!green] at (190,-450) {Battery};

% Wind 47
\coordinate[shift={(0,30)}] (wind) at (N47);
\draw[blue!50!green, very thick] (wind) -- (N47);
\draw[blue!50!green, fill=white, rotate=50] (wind)+(10,0) ellipse (10 and 1);
\draw[blue!50!green, fill=white, rotate=170] (wind)+(10,0) ellipse (10 and 1);
\draw[blue!50!green, fill=white, rotate=290] (wind)+(10,0) ellipse (10 and 1);
\draw[blue!50!green, fill=white] (wind) circle (3);

% Wind 48
\coordinate[shift={(0,30)}] (wind) at (N48);
\draw[blue!50!green, very thick] (wind) -- (N48);
\draw[blue!50!green, fill=white, rotate=50] (wind)+(10,0) ellipse (10 and 1);
\draw[blue!50!green, fill=white, rotate=170] (wind)+(10,0) ellipse (10 and 1);
\draw[blue!50!green, fill=white, rotate=290] (wind)+(10,0) ellipse (10 and 1);
\draw[blue!50!green, fill=white] (wind) circle (3);

% Wind 1
\coordinate[shift={(11,30)}] (wind) at (N1);
\draw[blue!50!green, very thick] (wind) -- (N1);
\draw[blue!50!green, fill=white, rotate=50] (wind)+(10,0) ellipse (10 and 1);
\draw[blue!50!green, fill=white, rotate=170] (wind)+(10,0) ellipse (10 and 1);
\draw[blue!50!green, fill=white, rotate=290] (wind)+(10,0) ellipse (10 and 1);
\draw[blue!50!green, fill=white] (wind) circle (3);

% Wind 4
\coordinate[shift={(-30,0)}] (wind) at (N4);
\draw[blue!50!green, very thick] (wind) -- (N4);
\draw[blue!50!green, fill=white, rotate=50] (wind)+(10,0) ellipse (10 and 1);
\draw[blue!50!green, fill=white, rotate=170] (wind)+(10,0) ellipse (10 and 1);
\draw[blue!50!green, fill=white, rotate=290] (wind)+(10,0) ellipse (10 and 1);
\draw[blue!50!green, fill=white] (wind) circle (3);

% Wind 6
\coordinate[shift={(0,30)}] (wind) at (N6);
\draw[blue!50!green, very thick] (wind) -- (N6);
\draw[blue!50!green, fill=white, rotate=50] (wind)+(10,0) ellipse (10 and 1);
\draw[blue!50!green, fill=white, rotate=170] (wind)+(10,0) ellipse (10 and 1);
\draw[blue!50!green, fill=white, rotate=290] (wind)+(10,0) ellipse (10 and 1);
\draw[blue!50!green, fill=white] (wind) circle (3);

% Wind 85
\coordinate[shift={(30,0)}] (wind) at (N85);
\draw[blue!50!green, very thick] (wind) -- (N85);
\draw[blue!50!green, fill=white, rotate=50] (wind)+(10,0) ellipse (10 and 1);
\draw[blue!50!green, fill=white, rotate=170] (wind)+(10,0) ellipse (10 and 1);
\draw[blue!50!green, fill=white, rotate=290] (wind)+(10,0) ellipse (10 and 1);
\draw[blue!50!green, fill=white] (wind) circle (3);

% Wind 87
\coordinate[shift={(20,30)}] (wind) at (N87);
\draw[blue!50!green, very thick] (wind) -- (N87);
\draw[blue!50!green, fill=white, rotate=50] (wind)+(10,0) ellipse (10 and 1);
\draw[blue!50!green, fill=white, rotate=170] (wind)+(10,0) ellipse (10 and 1);
\draw[blue!50!green, fill=white, rotate=290] (wind)+(10,0) ellipse (10 and 1);
\draw[blue!50!green, fill=white] (wind) circle (3);

% Wind 107
\coordinate[shift={(40,10)}] (wind) at (N107);
\draw[blue!50!green, very thick] (wind) -- (N107);
\draw[blue!50!green, fill=white, rotate=50] (wind)+(10,0) ellipse (10 and 1);
\draw[blue!50!green, fill=white, rotate=170] (wind)+(10,0) ellipse (10 and 1);
\draw[blue!50!green, fill=white, rotate=290] (wind)+(10,0) ellipse (10 and 1);
\draw[blue!50!green, fill=white] (wind) circle (3);

% Wind 109
\coordinate[shift={(-30,10)}] (wind) at (N109);
\draw[blue!50!green, very thick] (wind) -- (N109);
\draw[blue!50!green, fill=white, rotate=50] (wind)+(10,0) ellipse (10 and 1);
\draw[blue!50!green, fill=white, rotate=170] (wind)+(10,0) ellipse (10 and 1);
\draw[blue!50!green, fill=white, rotate=290] (wind)+(10,0) ellipse (10 and 1);
\draw[blue!50!green, fill=white] (wind) circle (3);

% Wind 113
\coordinate[shift={(0,30)}] (wind) at (N113);
\draw[blue!50!green, very thick] (wind) -- (N113);
\draw[blue!50!green, fill=white, rotate=50] (wind)+(10,0) ellipse (10 and 1);
\draw[blue!50!green, fill=white, rotate=170] (wind)+(10,0) ellipse (10 and 1);
\draw[blue!50!green, fill=white, rotate=290] (wind)+(10,0) ellipse (10 and 1);
\draw[blue!50!green, fill=white] (wind) circle (3);

% Leyenda Wind
\coordinate (wind) at (250,-430);
\draw[blue!50!green, very thick] (wind) -- +(0,-30);
\draw[blue!50!green, fill=white, rotate=50] (wind)+(10,0) ellipse (10 and 1);
\draw[blue!50!green, fill=white, rotate=170] (wind)+(10,0) ellipse (10 and 1);
\draw[blue!50!green, fill=white, rotate=290] (wind)+(10,0) ellipse (10 and 1);
\draw[blue!50!green, fill=white] (wind) circle (3);
\node[blue!50!green] at (280,-450) {Wind};

\node[blue!50!green] at (N117) [yshift=-8,blue!50!green]  {Slack};
\end{tikzpicture}

%% file: Figures/Results2_IEEE123.tex
\begin{tikzpicture}
\begin{axis}[width=9cm, height=5 cm, xlabel={Time},ylabel={Power (kW)},legend pos=north west, xmax=23, xmin=0, xmajorgrids,ymajorgrids,legend style={nodes={scale=0.6, transform shape},legend columns=4}]

\addplot[green!50!blue, thick, mark=*] coordinates {( 0 , 1035000.0 )
( 1 , 1035000.0 )
( 2 , 862500.0 )
( 3 , 862500.0 )
( 4 , 1035000.0 )
( 5 , 1207500.0 )
( 6 , 1552500.0 )
( 7 , 1897500.0 )
( 8 , 2242500.0 )
( 9 , 2415000.0 )
( 10 , 2587500.0 )
( 11 , 2760000.0 )
( 12 , 2587500.0 )
( 13 , 2587500.0 )
( 14 , 2587500.0 )
( 15 , 2587500.0 )
( 16 , 2415000.0 )
( 17 , 2415000.0 )
( 18 , 2932500.0 )
( 19 , 3450000.0 )
( 20 , 2932500.0 )
( 21 , 2242500.0 )
( 22 , 1897500.0 )
( 23 , 1207500.0 )};
\addlegendentry{Demand};

\addplot[blue!80!green, thick, mark=square] coordinates{( 0 , -0.014585488799639279 )
( 1 , -0.020675380072589178 )
( 2 , -0.010958808711620804 )
( 3 , -0.017457878548611916 )
( 4 , -0.012026192615621767 )
( 5 , -0.024551485802930983 )
( 6 , 116670.63788858801 )
( 7 , 721979.2506186239 )
( 8 , 1372582.1609819194 )
( 9 , 1192548.347091981 )
( 10 , 2102644.829392245 )
( 11 , 856577.7694161363 )
( 12 , 350246.324235338 )
( 13 , 766604.2243630973 )
( 14 , 597895.7309949007 )
( 15 , 697233.0135193061 )
( 16 , 637438.0997122765 )
( 17 , 942077.5906090941 )
( 18 , 1157873.259038546 )
( 19 , 1654627.404835723 )
( 20 , 1546721.9446686432 )
( 21 , 651533.3543410415 )
( 22 , 149246.12873197187 )
( 23 , -0.015498606672281312 )};
\addlegendentry{Power at PCC}

\end{axis}
\end{tikzpicture}

%% file: Figures/Results_IEEE123.tex
\begin{tikzpicture}
\begin{axis}[width=9cm, height=5 cm, xlabel={Time},ylabel={Power (kW)},legend pos=north west, xmax=23, xmin=0, xmajorgrids,ymajorgrids,legend style={nodes={scale=0.6, transform shape},legend columns=4}]

\addplot[orange, thick, mark=square] coordinates{( 0 , 0.0 )
( 1 , 0.0 )
( 2 , 0.0 )
( 3 , 0.0 )
( 4 , 0.0 )
( 5 , 0.0 )
( 6 , 104999.99991363155 )
( 7 , 209999.9999103565 )
( 8 , 314999.99991333287 )
( 9 , 419999.99990786746 )
( 10 , 498749.9999028222 )
( 11 , 524999.9998830883 )
( 12 , 472499.9999140818 )
( 13 , 446249.99990637356 )
( 14 , 419999.99990788614 )
( 15 , 314999.99991335993 )
( 16 , 209999.9999182941 )
( 17 , 104999.99992181436 )
( 18 , 0.0 )
( 19 , 0.0 )
( 20 , 0.0 )
( 21 , 0.0 )
( 22 , 0.0 )
( 23 , 0.0 )};
\addlegendentry{Solar};

\addplot[blue!50!green, thick, mark=*] coordinates{( 0 , 872260.3584249367 )
( 1 , 981054.2384831677 )
( 2 , 789068.8117843337 )
( 3 , 842917.9206228415 )
( 4 , 1025931.0779313594 )
( 5 , 1526268.2686604133 )
( 6 , 1148174.9998360756 )
( 7 , 1148174.9998361587 )
( 8 , 340199.99989190174 )
( 9 , 806399.9998667443 )
( 10 , 196874.99989810947 )
( 11 , 1148174.9998573014 )
( 12 , 1574999.9997959367 )
( 13 , 1574999.9997759508 )
( 14 , 1574999.9997756341 )
( 15 , 1574999.9997776072 )
( 16 , 1574999.9997797245 )
( 17 , 1574999.9997805017 )
( 18 , 1574999.9997931332 )
( 19 , 1574999.9997945386 )
( 20 , 1574999.999779667 )
( 21 , 1574999.9997789357 )
( 22 , 1574999.999780501 )
( 23 , 1420806.1686134767 )};
\addlegendentry{Wind};

\addplot coordinates{( 0 , 582741.3103097891 )
( 1 , 636687.541509682 )
( 2 , 710123.8456445385 )
( 3 , 729704.0936138054 )
( 4 , 738770.7213293555 )
( 5 , 419999.9540531315 )
( 6 , 602654.3598525485 )
( 7 , 419999.9537265637 )
( 8 , 634718.0560241125 )
( 9 , 630770.2961372964 )
( 10 , 419999.953842441 )
( 11 , 650246.7474608081 )
( 12 , 840000.0463874922 )
( 13 , 639644.7153829525 )
( 14 , 634248.7379394764 )
( 15 , 634515.0702428487 )
( 16 , 627076.173782422 )
( 17 , 419999.9535049258 )
( 18 , 619627.1091878101 )
( 19 , 840000.0465021462 )
( 20 , 650778.5481535129 )
( 21 , 666746.3123968834 )
( 22 , 840000.0469817875 )
( 23 , 626697.0954421947 )};
\addlegendentry{Battery};

\end{axis}
\end{tikzpicture}